\input amssym.def
\input amssym
\magnification=1100
\baselineskip = 0.25truein
\lineskiplimit = 0.01truein
\lineskip = 0.01truein
\vsize = 8.5truein
\voffset = 0.2truein
\parskip = 0.10truein
\parindent = 0.3truein
\settabs 12 \columns
\hsize = 5.4truein
\hoffset = 0.4truein

\setbox\strutbox=\hbox{%
\vrule height .708\baselineskip
depth .292\baselineskip
width 0pt}
\font\caps=cmcsc10

\def\sqr#1#2{{\vcenter{\vbox{\hrule height.#2pt
\hbox{\vrule width.#2pt height#1pt \kern#1pt
\vrule width.#2pt}
\hrule height.#2pt}}}}
\def\square{\mathchoice\sqr46\sqr46\sqr{3.1}6\sqr{2.3}4}

\def\PSL{{\rm PSL}}
\def\SL{{\rm SL}}
\def\PGL{{\rm PGL}}
\def\GL{{\rm GL}}
\def\P{{\rm P}}

\def\Ram{{\rm Ram}}
\def\rank{{\rm rank}}
\def\tr{{\rm tr}}

\centerline{\bf COVERING SPACES OF ARITHMETIC 3-ORBIFOLDS}
\vskip 6pt
\centerline{\bf {M. Lackenby, D. D. Long and A. W. Reid}}
\vskip 6pt
\centerline{\caps 1. Introduction}
\vskip 6pt
This paper investigates properties of finite
sheeted covering spaces of arithmetic hyperbolic 3-orbifolds (see \S 2).  
The main motivation is 
a central unresolved question in the theory of closed hyperbolic
3-manifolds; namely whether a closed hyperbolic 3-manifold is virtually
Haken.  Various strengthenings of this have also been widely studied.
Of specific to interest to us is the question of whether the fundamental group
of a given hyperbolic 3-manifold $M$ is {\sl large}; that is to say, some finite
index subgroup of $\pi_1(M)$
admits a surjective homomorphism onto a non-abelian free group. This implies that
$M$ is virtually Haken, and indeed that $M$ has {\sl infinite virtual
first Betti number} (see \S 2.4 for a definition). 
Of course, a weaker formulation is to only ask whether
the virtual first Betti number of
a closed hyperbolic 3-manifold $M$ is positive.
This has been verified in many cases, see [8] for some recent work
on this. However, in general, passage
from positive virtual first Betti number to infinite virtual first Betti number
is difficult, as is passage from infinite virtual first Betti number to
large.  This paper makes some progress on the latter in certain settings.

The background for our work is recent work of the first author
(see for example [17] and [19]). This suggests that the
questions addressed above 
for hyperbolic 3-manifolds that are commensurable with an
orbifold may be more amenable to study.  One of the aims of this paper is to
address these questions for arithmetic hyperbolic 3-manifolds 
and in particular, provide further evidence for a positive solution to
the largeness question.  It is already known that many arithmetic
hyperbolic 3-manifolds have infinite virtual first Betti number,
mainly through the application of the theory of automorphic forms (see [4],
[14], [20], [27] and [35]). For convenience, we shall refer to these 
collectively as
{\sl arithmetic methods}.  Some geometric methods are also known using
the existence of a totally geodesic surface, and largeness is known
there (see [21] and [25]).  However the question of largeness
remains unknown in general for arithmetic hyperbolic 3-manifolds, in 
particular even for those for which positive virtual first Betti number is
known by arithmetic methods.

Our main results are the following, the first of which explains why
arithmetic manifolds are particularly well-suited to the above
questions in the context of orbifolds.

\noindent {\bf Theorem 1.1.}~{\sl Let $M={\Bbb H}^3/\Gamma$ be an arithmetic hyperbolic 3-manifold. 
Then $M$ is commensurable 
with an arithmetic hyperbolic 3-orbifold with non-empty singular locus.  More
precisely, $\Gamma$ is commensurable with an arithmetic Kleinian group
$\Gamma_0$ containing an element of order $2$.}

In fact, more can be said about elements of order $2$ in the commensurability
class of $\Gamma$. 

\noindent {\bf Theorem 1.2.}~{\sl
Let $M={\Bbb H}^3/\Gamma$ be as above.
Then $\Gamma$ is commensurable with an arithmetic Kleinian group
$\Gamma_0$ containing a subgroup isomorphic to 
${\Bbb Z}/2{\Bbb Z}\times{\Bbb Z}/2{\Bbb Z}$. }

Theorem 1.1 is easily seen to
be false in the setting of non-arithmetic manifolds.  For example, by
Margulis's result [33], in the non-arithmetic case
there is a unique maximal
element in the commensurability class of the group. It is 
easy to construct examples whereby this maximal element 
has no non-trivial elements of finite order (see for example [36]).

We apply these results, together with results in \S 3,
 that are in the spirit of the [17] and [19]
to obtain the following results for arithmetic hyperbolic 3-manifolds.

\noindent {\bf Theorem 1.3.} {\sl Let $M$ be an arithmetic hyperbolic 3-manifold for which
the virtual first Betti number is at least $4$. Then $\pi_1(M)$ is large.}

In [1] Borel shows that arithmetic manifolds (not necessarily
hyperbolic) having a congruence subgroup with positive first Betti
number, have infinite virtual first Betti number (see \S 2.4). Thus we have.

\noindent {\bf Corollary 1.4.} {\sl Let $M$ be an arithmetic hyperbolic 3-manifold for which
arithmetic methods apply to produce a cover with positive first Betti number. Then $\pi_1(M)$ is large.}

In particular this applies to all known examples of arithmetic
hyperbolic 3-manifolds that have covers with positive first Betti
number (we discuss some specific examples of this in \S 6).

As further evidence for studying orbifolds, and in particular arithmetic ones,
we also show:

\noindent {\bf Theorem 1.5.} {\sl Let $M$ be an arithmetic hyperbolic 
3-manifold
commensurable with an orbifold $O={\Bbb H}^3/\Gamma$ 
such that either;
\item{(i)}$\Gamma$ contains $A_4$, $S_4$ or $A_5$ or;

\item{(ii)}$\Gamma$ is derived from a quaternion algebra and contains a finite
dihedral group.

Then $\pi_1(M)$ is large.}

The proof of Theorem 1.3 involves establishing linear growth in
mod $p$ homology for some prime $p$ (see \S 5) and has applications to
orbifolds other than arithmetic ones. 
For example we prove in \S 5 the following result.

\noindent {\bf Theorem 1.6.} {\sl Let $O$ be a 3-orbifold 
(with possibly empty singular locus) commensurable
with a closed orientable hyperbolic 3-orbifold that contains ${\Bbb Z}/2{\Bbb Z} \times {\Bbb Z}/2{\Bbb Z}$
in its fundamental group. Suppose that $vb_1(O) \geq 4$. Then $\pi_1(O)$
is large.}

Some of our other main results
concern this phenomena, and are independent of the results proved in \S 5.
For example in \S 4, we prove the following result. 

\noindent {\bf Theorem 1.7.} {\sl Let $O$ be a 3-orbifold 
(with possibly empty singular locus) commensurable with
a closed orientable hyperbolic 3-orbifold
that contains ${\Bbb Z}/2{\Bbb Z} \times {\Bbb Z}/2{\Bbb Z}$
in its fundamental group. Then $O$ has a tower of finite-sheeted
covers $\{ O_i \}$ that has linear growth of mod $2$ homology.}

This result is also proved in [17], using the Golod-Shafarevich
inequality and the theory of $p$-adic Lie groups.  Our proof uses only
properties of hyperbolic 3-orbifolds, and as such can be considered
``elementary''.  A consequence of this is that we can give a new
proof of the following result (see \S 10).  This was originally proved
by Lubotzky in [24], and again the proof used the Golod-Shafarevich
inequality, and the theory of $p$-adic Lie groups.

\noindent {\bf Theorem 1.8.} {\sl No arithmetic Kleinian group
has the congruence subgroup property.}

We also discuss Property $(\tau)$ in connection with orbifolds.
Property $(\tau)$ is an important group-theoretic
concept, introduced by Lubotzky and Zimmer [29]. It
has many applications to diverse areas of mathematics,
including hyperbolic 3-manifold theory (see [15] for more details). It
is conjectured that if $M$ is a closed orientable 3-manifold with
infinite fundamental group, then $\pi_1(M)$ does
not have Property $(\tau)$.
As is well-known, having virtually positive first Betti number implies
this, but beyond this, little is known by way establishing a group does
not have $(\tau)$. Another of our main results for arithmetic Kleinian groups
is.

\noindent {\bf Theorem 1.9.} {\sl Suppose that for every compact orientable 
3-manifold $M$ with infinite fundamental group, $\pi_1(M)$ 
fails to have Property $(\tau)$. Then any arithmetic Kleinian group
is large.}

\noindent{\bf Acknowledgements:}~{\it The first author is supported by the EPSRC
and the second and third authors are supported by the NSF. The third author
gratefully acknowledges the hospitality of Centre 
Interfacultaire Bernoulli of EPF Lausanne during the preparation of this work.}

\vskip 18pt
\centerline{\caps 2. Arithmetic hyperbolic 3-orbifolds}
\vskip 6pt

\noindent {\bf 2.1}~We begin by recalling some facts about arithmetic Kleinian
groups that will be needed (see [32] for further details).

Arithmetic Kleinian groups are obtained as follows. Let $k$ be a number 
field having exactly one complex place, and
$B$ a quaternion algebra
over $k$ which ramifies at all real places of $k$. 
Let $\rho :~B \rightarrow M(2,{\Bbb C})$ be
an embedding, $\cal O$ an order of $B$, and ${\cal O}^1$ the elements
of norm one in $\cal O$. Then $\P\rho({\cal O}^1) < \PSL(2,{\Bbb C})$ is a
finite co-volume Kleinian group, which is co-compact if and only if
$B$ is not isomorphic to $M(2,{\Bbb Q}(\sqrt{-d}))$, where $d$ is a square
free positive integer.
An {\sl arithmetic Kleinian group} $\Gamma$ is a subgroup of 
$\PSL(2,{\Bbb C})$ commensurable with
a group of the type $\P\rho({\cal O}^1)$. 
We call $Q = {\Bbb H}^3/\Gamma$ {\sl arithmetic}
if $\Gamma$ is arithmetic.

\noindent {\bf Notation:}~We shall 
denote $\P\rho({\cal O}^1)$ by $\Gamma_{\cal O}^1$ and the set of finite places of $k$ that ramify
the quaternion algebra $B$ by $\Ram_f(B)$.

An arithmetic Kleinian group 
$\Gamma$ is called {\sl derived from a quaternion algebra} if
$\Gamma < \Gamma_{\cal O}^1$. For convenience we
state the following result that is deduced from the characterisation theorem
for arithmetic Kleinian groups (see [32] Corollary 8.3.5). For a finitely
generated group $G$ we denote by $G^{(2)}$ the subgroup
of $G$ generated by the squares of elements in $G$.

\noindent{\bf Theorem 2.1.}~{\sl
Let $\Gamma$ be a finite co-volume Kleinian group. Then $\Gamma$ 
is arithmetic if and only if the group $\Gamma^{(2)}$ is derived
from a quaternion algebra.}

One final fact about arithmetic Kleinian groups that we will make use
of is the following. If $\Gamma$ is derived from a quaternion algebra $B$
defined over $k$, then 
$$A_0\Gamma=\{\sum a_i \gamma_i~:~a_i \in k, \gamma_i \in \Gamma\},$$
is a quaternion algebra over $k$ (see [32] Chapter 3) and
is isomorphic to $B$ (see [32] Chapter 8).
In what follows we shall just identify the two.

\noindent{\bf Remark:}~For convenience, we have blurred the distinction
between an element $a\in \PSL(2,{\Bbb C})$ and a matrix $A\in \SL(2,{\Bbb C})$
that projects to $a$ under the homomorphism 
$\SL(2,{\Bbb C})\rightarrow \PSL(2,{\Bbb C})$. 
        
\noindent {\bf 2.2}~ Here we prove Theorems 1.1 and 1.2 (Theorem 1.1 is
implicit in [11]).

We begin with a lemma.

\noindent {\bf Lemma 2.2.}~{\sl
Let $\Gamma$ be derived from a quaternion algebra $B$, defined over
the number field $k$. Let $R_k$ denote the ring of integers of $k$,
$a$ and $b$ a pair of non-commuting elements of $\Gamma$, and
let ${\cal O} = R_k[1,a,b,a.b]$. Then $\cal O$ is an order of $B$.}

\noindent {\sl Proof.}~To show that $\cal O$ is an order we proceed
as follows.  First, since $a$ and $b$ do not commute, it is easy to
see that $\{1,a,b,ab\}$ spans $B$ over $k$. Thus $R_k[1,a,b,ab]$ contains a
$k$-basis of $B$, is finitely generated and contains $R_k$. Also note that
since $\Gamma$ is derived from a quaternion algebra the elements $a$,
$b$, $a.b$ and $R_k$-combinations of these words are integral in the
algebra.  To complete the proof, it suffices to prove that all
products of the basis elements can be expressed as $R_k$-combinations
of the basis elements. This follows from the Cayley-Hamilton theorem
as well as some other trace identities that we include
below.
\smallskip
\noindent $a+a^{-1} = {\tr }(a)1$,
\smallskip
\noindent $a^2  =  {\tr} (a)a - 1$, 
\smallskip
\noindent $a^{2} b  =  {\tr} (a)ab - b$, 
\smallskip
\noindent $aba  = -{\tr} (b)1+{\tr} (ab)a+b$, 
\smallskip
\noindent $b^{-1}a^{-1} = {\tr}(b)a^{-1}-ba^{-1}$,
\smallskip
\noindent $ba+ab  = ({\tr} (ab)-{\tr}(a){\tr}(b))1+{\tr}(b)a
                                    +{\tr}(a)b$.
\smallskip
In particular note that the first identity, with $a$ replaced by $ab$ throughout,
implies that $b^{-1}a^{-1} \in {\cal O}$ and the last identity then implies
that $a^{-1}b^{-1}\in {\cal O}$. $\square$

\noindent{\bf Remark.}~The discriminant of the order $\cal O$ in Lemma 2.2
can be easily computed and is the ideal $<\tr[a,b]-2>$.

Define the {\sl normalizer} of $\cal O$ in $B$ by:
$$N({\cal O}) = \{x \in B^* \mid x{\cal O}x^{-1} = {\cal O}\}.$$

The image, $\Gamma({\cal O})$ of $N({\cal O})$ in $\PGL(2,{\Bbb C})$
(which is isomorphic to $\PSL(2,{\Bbb C})$), is an arithmetic Kleinian
group. To see this we argue as follows.

Note first that, for every $x\in {\cal O}^1$, 
$x {\cal O} x^{-1} = {\cal O}$ because ${\cal O}$ is a
ring. Hence $N({\cal O})$ contains ${\cal O}^1$.  
Furthermore, any element of $N({\cal O})$
normalizes ${\cal O}^1$ (since conjugation preserves the norm).
Therefore $\Gamma({\cal O})$ is a subgroup of the normalizer of
$\Gamma_{\cal O}^1$ in $\PGL(2,{\Bbb C})$. Since $\Gamma_{\cal O}^1$
has finite co-volume, it is well-known that its normalizer is also
discrete and finite co-volume. Hence, $\Gamma({\cal O})$ is discrete. 
It also has finite co-volume, since it contains $\Gamma_{\cal O}^1$.
We summarize this discussion in the
following.

\noindent {\bf Corollary 2.3.}~{\sl
Let $B$ and $\cal O$ be as above. Then $\Gamma({\cal O})$ is an arithmetic
Kleinian group commensurable with $\Gamma_{\cal O}^1$.}

\noindent{\bf Remark.}~Corollary 2.3 holds more
generally. Namely, if $\cal O$ is any order of a quaternion
algebra $B$ (as in \S 2.1), then $N({\cal O})$ always gives rise to an
arithmetic Kleinian group $\Gamma({\cal O})$ (see [32] Chapter
6). We have included the above proof for completeness, and since it is 
straightforward in this case.

Theorem 1.1 will follow immediately from
the next proposition and Theorem 2.1. This will require some
notation.

Let $a$ and $b$ be elements of $\SL(2,{\Bbb C})$ without a common fixed point.  
Then, as noticed by J\o rgenson [12], $ab-ba$ is 
an element of $\GL(2,{\Bbb C})$ which has trace 0 and whose image
in $\PGL(2,{\Bbb C})$ is of order two and conjugates $a$ to $a^{-1}$
and $b$ to $b^{-1}$.  Denote this involution by $\tau_{a,b}$.

\noindent{\bf Proposition 2.4.}~{\sl
Let $\Gamma$ be derived from a quaternion algebra and $a,b\in\Gamma$
such that $H = < a,b >$ is a non-elementary subgroup of $\Gamma$.
Then $\tau_{a,b}$ is contained in an arithmetic
Kleinian group commensurable with $\Gamma$.}

\noindent {\sl Proof.}~Since $\Gamma$ is derived from a quaternion
algebra, there exists an order $\cal D$ (as in \S 2.1) such that
$\Gamma < \Gamma_{\cal D}^1$. Let ${\cal O} = R_{k\Gamma}[1, a, b, a.b]$ 
be as in Lemma 2.2. Note that $\cal O \subset \cal D$.
By Corollary 2.3, $\Gamma({\cal O})$
is an arithmetic Kleinian group that is commensurable with $\Gamma_{\cal O}^1$.
This in turn is commensurable with $\Gamma_{\cal D}^1$, and hence $\Gamma$.
Finally, the involution $\tau_{a,b}\in \Gamma({\cal O})$. To see this 
note that $\tau_{a,b}(ab) = a^{-1}b^{-1}$ which is an element of $\cal O$
by Lemma 2.2. $\square$

\noindent{\sl Proof of Theorem 1.2.}
The extension of the argument to prove Theorem 1.2 is made as
follows. $\Gamma$ will continue to be derived from a quaternion
algebra $B$ and we choose $a$ and $b$ loxodromic elements
such that their axes, $A_a$ and $A_b$ respectively,
are disjoint. 

By construction, the involution
$\tau_{a,b}$ rotates around the geodesic $\gamma_{a,b}$
that is the common perpendicular between $A_a$ and $A_b$.  
We now claim that there is an involution $\tau_{\alpha_,\beta}$ that
acts by rotating around $A_a$.  To prove the claim, first observe that
since $\Gamma({\cal O})$ (as in the proof of Theorem 1.1) has
finite co-volume, there is a loxdromic element in $\Gamma({\cal O})$
that has $\gamma_{a,b}$ as an axis.  Since $\Gamma$ is commensurable
with $\Gamma({\cal O})$ there is a loxodromic element $\alpha\in
\Gamma$ that has $\gamma_{a,b}$ as an axis. Note that since 
$\alpha$ and $\tau_{a,b}$ share an axis, they commute.
In addition, there is a
loxodromic element $\beta\in \Gamma$ that has the geodesic
$a\gamma_{a,b}$ as an axis.  Hence, $A_a$ is the common perpendicular
of the axes $\gamma_{a,b}$ and $a\gamma_{a,b}$, and as before we can
construct an involution $\tau_{\alpha,\beta}$ that acts as claimed (and
commutes with $a$).

Note that $\tau_{a,b}$ and $\tau_{\alpha,\beta}$ are involutions and
commute. Hence,
the group $V=<\tau_{a,b},\tau_{\alpha,\beta}>$ is isomorphic to 
${\Bbb Z}/2{\Bbb Z}\times{\Bbb Z}/2{\Bbb Z}$. It remains to show that $V$
is a subgroup of an arithmetic Kleinian group commensurable with $\Gamma$.
To see this let $\cal L$ be the order associated to the group $<a,\alpha>$
as in Lemma 2.2. The action of these involutions on $a$ and $\alpha$ is
given by:

$$\tau_{a,b}a\tau_{a,b} = a^{-1},~~~\tau_{a,b}\alpha\tau_{a,b} = \alpha,$$
$$\tau_{\alpha,\beta}a\tau_{\alpha,\beta} = a,~~~
           \tau_{\alpha,\beta}\alpha\tau_{\alpha,\beta} = \alpha^{-1}.$$
It follows from this that $V<\Gamma({\cal L})$ and Corollary 2.3
completes the proof. $\square$

\noindent {\bf 2.3}~In this subsection we discuss implications on the
Hilbert symbol of the invariant quaternion algebra associated to a
Kleinian group of finite co-volume given the presence of $A_4$, $S_4$
or $A_5$ subgroup and certain dihedral subgroups. In the case of $S_4$
and $A_5$, since both of these contain $A_4$, we will restrict
consideration to this group.

\noindent{\bf Definition.}~Let $G$ be a finite subgroup of an arithmetic
Kleinian group.  We shall call $G$ {\sl derived from a quaternion
algebra} if $G$ is contained in some group $\Gamma_{\cal O}^1$ 
as above.

\noindent{\bf Theorem 2.5}~{\sl Suppose that $\Gamma$ is an arithmetic Kleinian
group commensurable with a Kleinian group containing $A_4$ or a Kleinian
group containing a finite dihedral group derived from a quaternion algebra.
Let $k$ and $B$ denote the
invariant trace-field and quaternion algebra of $\Gamma$. Then 
if $\nu\in \Ram_f B$ and $\nu$ divides the rational prime $p$, then 
$k_\nu$ contains no quadratic extension of ${\Bbb Q}_p$.}

\noindent{\bf Proof.}~Note first that if $\Gamma$ is commensurable with
a group $\Gamma_1$ containing $A_4$, then since $A_4=A_4^{(2)}$
it follows that $A_4<\Gamma_1^{(2)}$ and so any $A_4$ is derived 
from a quaternion algebra.  Furthermore, $A_4$ contains a copy
of ${\Bbb Z}/2{\Bbb Z}\times{\Bbb Z}/2{\Bbb Z}$, which is the dihedral group
of order $4$. Thus we can assume that we are in the case that $\Gamma$
is a Kleinian group derived from a quaternion algebra and contains a dihedral
group $D_{n}$ of order $2n$.

We can assume that $\Gamma$ is cocompact, otherwise, $B$ is a matrix algebra
and is unramified at all places of $k$.
Let $x,y\in \Gamma$ generate the dihedral subgroup, with

$$D_n = <x,y|x^2=y^2=(xy)^n=1>.$$

Note that since $x$ and $y$ do not have a common fixed point on the sphere at
infinity, it follows that
a Hilbert symbol for $B$ can be computed using the
basis $\{1,x,y,xy\}$. From [32] Theorem 3.6.1 we deduce that a Hilbert
symbol is given by

$$B \cong \biggl({{-4, 4\cos^2 2\pi/n - 4}\over k}\biggr)
 \cong \biggl({{-1, 4\cos^2 2\pi/n - 4}\over k}\biggr)$$

We need the following information about
the term $\tau_n=4\cos^2 2\pi/n - 4$ (cf. [31] Lemma 4.4).

\noindent {\bf Lemma 2.6.} {\sl If $n$ is odd, or even and greater than $4$,
then $\tau_n$ has norm $p$ or is a unit,
depending on whether $n$ is a power of a prime $p$ or not.}

In the case $n=4$, $\tau_n=-4$, and so the Hilbert symbol becomes
$\biggl({{-1, -4}\over k}\biggr)\cong \biggl({{-1, -1}\over k}\biggr)$.

Given this, and the lemma, we gain some preliminary control on $\Ram_f B$.
For, if $\nu\in \Ram_f B$ then from above 
we deduce that 
$\nu$ divides $2$ or at most one other rational prime $p$
(see [32] Theorem 2.6.6). 
Furthermore, the order ${\cal O} = R_k[1,x,y,xy]$ (recall Lemma 2.2)
can be shown to have discriminant $d({\cal O}) = <\tau_n>$ (see the
Remark following Lemma 2.2). Now, the discriminant of a maximal order of $B$, 
which equals the product of finite places ramifying $B$, divides
$d({\cal O}) = <\tau_n>$ (see [32] Theorem 6.3.4). Thus it follows
that if $p$ is odd, then $\nu$ cannot divide 2. Given these remarks
we can now argue as follows.

\medskip

\noindent{\bf Case 1:}~Assume $n$ is not a prime power, and so $\tau_n$
is a unit. Hence $d({\cal O})$ is the trivial ideal, and so
it follows from the discussion
above that $\cal O$ is maximal. Hence $B$ is unramified at all finite
places, and the theorem is proved in this case.

\medskip

\noindent{\bf Case 2:}~Assume $n=p^t$ is a prime power, and $n\neq 4$. 
>From the remarks
preceeding Case 1, it follows that $\Ram_f B =\emptyset$ or consists of 
a unique place 
dividing $p$. In particular
if $p\neq 2$, it cannot contain places dividing $2$.
Assume
by way of contradiction, that $\Ram_f B$ contains a place $\nu$
such that $k_\nu$ contains a quadratic extension $\ell$ of ${\Bbb Q}_p$.

Now ${\Bbb Q}(\cos 2\pi/{p^t})$ is a subfield of $k$, and so
$B$ can be described as follows.

$$B \cong \biggl({{-1, \tau_{p^t}}\over{\Bbb Q}(\cos 2\pi/{p^t})}\biggr)
                               \otimes_{{\Bbb Q}(\cos 2\pi/{p^t})} k.$$

Assume first that $[k:{\Bbb Q}]$ has even degree.
Hence the quaternion algebra $B$ is ramified at all real places of $k$ 
(an even number). Hence if $\Ram_f B$ is non-empty it consists of an even
number of finite places. However, recall
from above that the order has discriminant $d({\cal O})$ which is 
either the trivial ideal or a prime ideal of norm $p$. 
As this discriminant
divides that of a maximal order, and the cardinality of $\Ram_f B$ is even,
it follows that the discriminant of a maximal order,
and hence $B$, is the trivial ideal. Hence we are done in this case.

Now assume that $[k:{\Bbb Q}]$ is odd, and
so $[{\Bbb Q}(\cos 2\pi/{p^t}):{\Bbb Q}]$ is odd.
By the theory of ramification of
primes in the maximal real subfield of a cyclotomic field, there is a unique 
${\Bbb Q}(\cos 2\pi/{p^t})$-prime $\omega$ dividing $p$ and this has norm $p$.
Furthermore, since ${\Bbb Q}(\cos 2\pi/{p^t})$ is assumed to have 
odd degree, the theory of ramification in number field extensions [34] 
implies that ${\Bbb Q}(\cos 2\pi/{p^t})_\omega$  also will
have odd degree, and so must be disjoint from the field $\ell$.
Hence $k_\nu$ contains a subfield $L$ that is the compositum of 
${\Bbb Q}(\cos 2\pi/{p^t})$ and $\ell$.
$L$ has degree 2 over ${\Bbb Q}(\cos 2\pi/{p^t})$ 
and so applying [32] Theorem 2.6.5, $L$ must split the algebra 
$\biggl({{-1, \tau_{p^t}}\over{\Bbb Q}(\cos 2\pi/{p^t})_\omega}\biggr)$.  
However, then we have from above,
$$B_\nu  = B\otimes_k k_\nu \cong 
\biggl({{-1, \tau_{p^t}}\over{\Bbb Q}(\cos 2\pi/{p^t})_\omega}\biggr) 
\otimes_{{\Bbb Q}(\cos 2\pi/{p^t})_\omega} L \otimes_L k_\nu \cong M(2,k_\nu)$$
which is a contradiction, since $\nu$ ramifies.

\noindent {\bf Case 3:} Finally, we must deal with the case $n=4$. This is similar to the case above.
In this case, if $\nu\in \Ram_f(B)$ then $\nu$ is dyadic,
and so $k_\nu$ is a finite extension of ${\Bbb Q}_2$.
Suppose there is  a quadratic 
extension $\ell$ of ${\Bbb Q}_2$ that is contained in $k_\nu$.
Now $B$ is ramified at $\nu$, and so $B_\nu = B\otimes_k k_\nu$ is isomorphic to the unique
division algebra over $k_\nu$. As above,
it is easy to see that the following tensor
products hold:
$$B_\nu = B\otimes_k k_\nu = \biggl ({{-1,-1}\over k_\nu}\biggr)
\cong \biggl({{-1,-1}\over {{\Bbb Q}_2}}\biggr) \otimes_{{\Bbb Q}_2} k_\nu.$$
$$\cong \biggl({{-1,-1}\over {{\Bbb Q}_2}}\biggr) \otimes_{{\Bbb Q}_2} \ell
                          \otimes_\ell k_\nu$$

However, arguing as above,
a quadratic extension of the center of a quaternion algebra over a local field
(in particular $\ell$ over ${\Bbb Q}_2$) splits the unique
division algebra over ${\Bbb Q}_2$. Hence,

$$ \biggl({{-1,-1}\over {{\Bbb Q}_2}}\biggr) \otimes_{{\Bbb Q}_2}
             \ell \cong M(2,\ell).$$

$$B_\nu = B\otimes_k k_\nu \cong M(2,\ell) \otimes_\ell k_\nu
                         \cong M(2,k_\nu),$$
Again, this contradicts the assumption that $B_\nu$ is a division algebra. $\square$

\noindent {\bf 2.4}~In this subsection we gather together some notions
and results pertaining to congruence covers of arithmetic hyperbolic 
3-orbifolds. 
If $\Gamma$ is an arithmetic Kleinian group, 
there is a distinguished class of subgroups in $\Gamma$, known
as the {\sl congruence subgroups}. These
are defined as follows. Notation as in \S 2.1.

Let $\cal O$ be a maximal order of $B$, and 
let $I$ be any proper 2-sided integral ideal of $B$ contained in 
${\cal O}$; ie
$I$ is a complete $R_k$-lattice in $B$ such that
$${\cal O} = \{ x \in B \mid x I \subset I \} = \{ x \in B \mid I x \subset I \} .$$
As noted in [32] Chapter 6.1,
any proper 2-sided integral ideal of $B$ contained in 
${\cal O}$ is an ideal of $\cal O$ in the usual non-commutative ring sense.
In particular ${\cal O}/I$ is a finite ring.

Define
$${\cal O}^1(I) = \{ \alpha \in {\cal O}^1 : \alpha - 1 \in I \}.$$
The corresponding {\sl principal congruence subgroup} 
of $\Gamma_{\cal O}^1$ is 
$$\Gamma({\cal O}(I)) = \P\rho({\cal O}^1(I)).$$

If $\Gamma$ is an arithmetic Kleinian group then a subgroup
$\Delta < \Gamma$ is a {\sl congruence subgroup}
of $\Gamma$ if it contains some principal congruence 
subgroup $\Gamma({\cal O}(I))$ as above.

Before stating the result about the first Betti number of congruence subgroups
we require, we need some notation.

\noindent {\bf Notation.}~If $X$ is a group, space or orbifold, we
will denote by $b_1(X)$ the rank of $H_1(X;{\Bbb Z})\otimes
{\Bbb Q}$ and set
$$vb_1(X) = \sup \{ b_1(\tilde X) : \tilde X \hbox{ is a finite
index subgroup or finite cover of } X \}.$$

\noindent{\bf Theorem 2.7.} (Borel [1])~{\sl Suppose $\Gamma$ is an
arithmetic Kleinian group containing a subgroup
$\Gamma({\cal O}(I))$ with $b_1(\Gamma({\cal O}(I))>0$.
Then $vb_1(\Gamma)=\infty$.}

\vskip 18pt

\centerline{\caps 3. The homology of 3-orbifolds}
\vskip 6pt

\noindent {\bf Definition.} Let $O$ be a compact
orientable 3-orbifold.
Let ${\rm sing}(O)$ be its singular locus,
and let $|O|$ denote its underlying 3-manifold.
Let ${\rm sing}^0(O)$ and ${\rm sing}^-(O)$
denote the components of the
singular locus with, respectively, 
zero and negative Euler characteristic.
For any prime $p$,
let ${\rm sing}_p(O)$ denote the union of the
arcs and circles in ${\rm sing}(O)$
with singularity order that is a multiple of $p$.
Let ${\rm sing}_p^0(O)$ and ${\rm sing}_p^-(O)$
denote those components of ${\rm sing}_p(O)$
with zero and negative Euler characteristic.

When $O$ is closed, ${\rm sing}(O)$ is a disjoint union
of simple closed curves and trivalent graphs,
and hence ${\rm sing}(O) = {\rm sing}^0(O)
\cup {\rm sing}^-(O)$. However, it need not be the case
that ${\rm sing}_p(O) = {\rm sing}_p^0(O)
\cup {\rm sing}_p^-(O)$.

\noindent {\bf Terminology.} If $p$ is a prime, let
${\Bbb F}_p$ denote the field of order $p$.
If $X$ is a group, space or orbifold, let $d_p(X)$
be the dimension of $H_1(X; {\Bbb F}_p)$.

The following lower bound on homology will be a crucial tool
that we use throughout the rest of this paper.

\noindent {\bf Proposition 3.1.} {\sl Let $O$ be a
compact orientable 3-orbifold, and let $p$ be
a prime. Then $d_p(O) \geq b_1({\rm sing}_p(O))$.}

\noindent {\sl Proof.} Let $M$ denote the
3-manifold obtained from $O$ by a removing
an open regular neighbourhood of its singular
locus. Let $\{ \mu_1, \dots, \mu_r \}$ be a collection
of meridian curves, one encircling each
arc or circle of the singular locus.
Let $n_i$ be the singularity order of the arc
or circle that $\mu_i$ encircles. Then 
$$\pi_1(O) = \pi_1(M)/ \langle \! \langle
\mu_1^{n_1}, \dots, \mu_r^{n_r} \rangle \! \rangle.$$
Hence, 
$$H_1(O; {\Bbb F}_p) = H_1(M; {\Bbb F}_p) 
/ \langle \! \langle
\mu_1^{n_1}, \dots, \mu_r^{n_r} \rangle \! \rangle.$$
Now, when $n_i$ is coprime to $p$, quotienting
$H_1(M; {\Bbb F}_p)$
by $\mu_i^{n_i}$ is the same as quotienting by
$\mu_i$. And when $n_i$ is a multiple of
$p$, then quotienting $H_1(M; {\Bbb F}_p)$ by $\mu_i^{n_i}$ 
has no effect. Thus, if we let $M'$ be the
3-manifold obtained from $|O|$ by removing
an open regular neighbourhood of ${\rm sing}_p(O)$,
then $d_p(O) = d_p(M')$.
Now, it is a well known consequence of Poincar\'e
duality that, for the compact orientable 3-manifold $M'$,
$d_p(M') \geq {1 \over 2} d_p(\partial M') \geq
b_1({\rm sing}_p(O))$, as required. $\square$

\vskip 18pt
\centerline{\caps 4. Linear growth of homology}
\vskip 6pt

\noindent {\bf Definition.} Let $X$ be a group, space or orbifold and let $p$ be a prime.
Then a collection $\{ X_i \}$ of finite index subgroups
or finite-sheeted covers of $X$ with index or degree $[X:X_i]$
is said to have {\sl linear growth of mod $p$ homology}
if
$$\inf_i d_p(X_i)/[X:X_i] > 0.$$

In this section we prove Theorem 1.7, which we restate below for convenience.

\noindent {\bf Theorem 4.1.} {\sl Let $O$ be a 
3-orbifold (with possibly empty singular locus) commensurable with a 
closed orientable hyperbolic
3-orbifold that contains ${\Bbb Z}/2{\Bbb Z} \times {\Bbb Z}/2{\Bbb Z}$
in its fundamental group. Then $O$ has a tower of finite-sheeted
covers $\{ O_i \}$ that has linear growth of mod $2$ homology.}

This is a consequence of the following more general theorem
(Theorem 1.1 in [17]).

\noindent {\bf Theorem 4.2.} {\sl Let $O$ be a compact orientable
3-orbifold with non-empty singular locus and a finite-volume
hyperbolic structure. Let $p$ be a prime that divides the
order of an element of $\pi_1(O)$. Then $O$
has a tower of finite-sheeted covers $\{ O_i \}$ that
has linear growth of mod $p$ homology.}

However, the proof of Theorem 4.2 required some results about
$p$-adic analytic groups and it also used the Golod-Shafarevich
inequality. In this section, we will provide a much simpler
proof of the weaker Theorem 4.1.

Note first that we can assume that $O$ itself is closed, orientable and
hyperbolic and contains
${\Bbb Z}/2{\Bbb Z} \times {\Bbb Z}/2{\Bbb Z}$ in its fundamental group.
For, we know that $O$ is commensurable
with some such orbifold $O'$. Let $O''$ be a common cover
of $O$ and $O'$. We may assume that $O''$ is a regular cover of $O'$.
Suppose that we could prove Theorem 4.1 for $O'$,
providing a sequence of covers $\{ O_i \}$ with linear
growth of mod 2 homology. Then the covering spaces of $O'$
corresponding to $\pi_1(O'') \cap \pi_1(O_i)$ also have linear
growth of mod 2 homology, by the following elementary result, which appears as
Lemma 5.3 in [17]. 

\noindent {\bf Lemma 4.3.} {\sl Let $\{ G_i \}$ be a sequence of
finite index subgroups of a finitely generated group $G$,
and let $H$ be a finite index normal subgroup of $G$. If $\{ G_i \}$ has linear
growth of mod $p$ homology for some prime $p$, then $\{ G_i \cap H \}$
does also.}

So, let us suppose that $O$ is closed, orientable and
hyperbolic and contains
${\Bbb Z}/2{\Bbb Z} \times {\Bbb Z}/2{\Bbb Z}$ in its fundamental group.
We will prove Theorem 4.1 by finding a tower of finite
covers $O_i$ such that
$$\inf_i b_1({\rm sing}_2^-(O_i))/[O:O_i] >0.$$
By Proposition 3.1, $d_2(O_i)$ is at least $b_1({\rm sing}_2(O_i))$,
which is, of course, at least $b_1({\rm sing}_2^-(O_i))$.
Thus, $\{ O_i \}$ will indeed have linear growth of mod
2 homology.

The first step is to find a finite cover 
$\tilde O$ of $O$ such that ${\rm sing}_2^-(\tilde O)$
is non-empty.

\noindent {\bf Proposition 4.4.} {\sl Let $O$ be a closed
orientable hyperbolic 3-orbifold that contains ${\Bbb Z}/2{\Bbb Z} \times {\Bbb Z}/2{\Bbb Z}$
in its fundamental group. Then $O$ is finitely covered
by a 3-orbifold $\tilde O$ such that every arc and circle of
${\rm sing}(\tilde O)$ has order 2, and which contains
at least one singular vertex. In particular,
${\rm sing}_2^-(\tilde O)$ is non-empty.}

\noindent {\sl Proof.} Since $O$ is hyperbolic, Selberg's lemma
implies that it has a finite-sheeted regular cover that is a manifold
$M$. Let $\tilde O$ be the cover of $O$ corresponding to the subgroup
$\pi_1(M) ( {\Bbb Z}/2{\Bbb Z} \times {\Bbb Z}/2{\Bbb Z})$ of $\pi_1(O)$.
Then $M$ regularly covers $\tilde O$ with covering group
$\pi_1(\tilde O)/\pi_1(M) = ( {\Bbb Z}/2{\Bbb Z} \times {\Bbb Z}/2{\Bbb Z}) / 
(\pi_1(M) \cap ( {\Bbb Z}/2{\Bbb Z} \times {\Bbb Z}/2{\Bbb Z})) = {\Bbb Z}/2{\Bbb Z} \times {\Bbb Z}/2{\Bbb Z}$. 
Thus, any arc or circle of the singular locus of $\tilde O$ has order 2.
Hence, ${\rm sing}_2(\tilde O) = {\rm sing}(\tilde O)$.
Since $O$ is closed, ${\rm sing}(O)$ consists of
simple closed curves and trivalent graphs.
Now, $\tilde O$ is hyperbolic and so is obtained as the quotient
of ${\Bbb H}^3$ by the action of $\pi_1(\tilde O)$. The
${\Bbb Z}/2{\Bbb Z} \times {\Bbb Z}/2{\Bbb Z}$ subgroup of $\pi_1(\tilde O)$
is realised by a finite subgroup of ${\rm Isom}({\Bbb H}^3)$,
which must have a common fixed point in ${\Bbb H}^3$.
The image of this point in $\tilde O$ is a singular
vertex. Hence, some component of ${\rm sing}_2(\tilde O)$ therefore
has negative Euler characteristic. $\square$

The next step is to pass to a finite cover $O_1$
such that $b_1({\rm sing}_2^-(O_1))$ is arbitrarily
large. Note that, for any finite cover $O_1$ of $\tilde O$,
${\rm sing}(O_1) = {\rm sing}_2(O_1)$. Thus, ${\rm sing}_2(O_1)$
consists of trivalent graphs and simple closed curves. So, 
$$b_1({\rm sing}_2^-(O_1)) \geq |V({\rm sing}_2(O_1))|/2 + 1,$$
where $V({\rm sing}_2(O_1))$ is the vertices of the singular set.
Thus, we will establish a lower bound on $b_1({\rm sing}_2^-(O_1))$
by finding a lower bound on the number of singular vertices
of $O_1$.

\noindent {\bf Theorem 4.5.} {\sl Let $O$ be a closed orientable
hyperbolic 3-orbifold that contains 
${\Bbb Z}/2{\Bbb Z} \times {\Bbb Z}/2{\Bbb Z}$ in its fundamental group. Then, for any integer $N$,
$O$ has a finite-sheeted cover $O_1$ such that 
each arc and circle of ${\rm sing}(O_1)$ has order 2,
and which contains at least $N$ singular vertices. Hence, 
$b_1({\rm sing}_2^-(O_1)) \geq N/2+1$.}

The following is a key step in the proof of this theorem.

\noindent {\bf Theorem 4.6.} {\sl Let $O$ be a compact orientable
hyperbolic 3-orbifold (with possibly empty singular locus), 
and let $n$ be a positive integer. Then for infinitely many 
$n$-tuples of distinct primes $(p_1, \dots, p_n)$, $\pi_1(O)$ admits a 
surjective homomorphism $\phi$ onto $\prod_{i=1}^n \PSL(2,p_i)$.
Furthermore, if $\pi_i \colon \prod_{i=1}^n \PSL(2,p_i) \rightarrow
\PSL(2, p_i)$ is projection onto the $i^{\rm th}$ factor, then
we may ensure that ${\rm ker}(\pi_i \phi)$ is torsion-free,
for each $i$.}

\noindent{\sl Proof.} Let $O= {\Bbb H}^3/\Gamma$. It is shown in
[23] that for infinitely many rational primes $p$ there are (reduction)
homomorphisms $\phi_p \colon \Gamma \rightarrow \PSL(2,p)$.  
It is well-known
that by avoiding a finite set of primes we can assume that the kernels
are torsion-free (see Lemma 6.5.6 of [32] for example).
Also, by
definition of these homomorphisms, for all non-trivial elements $g \in \Gamma$,
$\phi_p(g) \neq 1$ for all but a finite number of primes.
Let $J$ be the set of rational primes $p$ given by the above construction.
It also follows from the argument in [23] that, for any finite index subgroup of $\Gamma$,
the restriction of $\phi_p$ to that subgroup is a surjection onto 
$\PSL(2,p)$, for all but finitely many primes $p$ in $J$. The proof 
is completed using a result of P. Hall [9] which asserts:

If $\Gamma$ is a group
and $\phi_i \colon \Gamma \rightarrow G_i$ are epimorphisms 
to distinct non-abelian finite
simple groups, then the product homomorphism 
$\Gamma \rightarrow \prod G_i$ is onto.$\square$

\noindent {\sl Proof of Theorem 4.5.} Let $O = {\Bbb H}^3/\Gamma$.
By Proposition 4.4, we may assume that
the order of each arc and circle of ${\rm sing}(O)$ is 2,
and that ${\rm sing}(O)$ contains at least one vertex. Let $\phi \colon 
\Gamma \rightarrow \prod_{i=1}^n {\rm PSL}(2, p_i)$
be the homomorphism from Theorem 4.6, 
let $H$ be the image of
${\Bbb Z}/2{\Bbb Z} \times {\Bbb Z}/2{\Bbb Z}$ under $\phi$, and
let $(A_1, \dots, A_n)$ and $(B_1, \dots, B_n)$ be the images under 
$\phi$ of the generators of ${\Bbb Z}/2{\Bbb Z} \times {\Bbb Z}/2{\Bbb Z}$. 
Since ${\rm ker}(\pi_i \phi)$ is torsion-free, $A_i$, $B_i$ and
$A_i B_i$ are non-trivial for each $i$. Let $O_1={\Bbb H}^3/\Gamma_1$ be 
the covering space of $O$ corresponding to the subgroup
$\phi^{-1}(H) = \Gamma_1$, and let $M$ be the covering space corresponding
to the kernel of $\phi$. Since the kernel of $\phi$ is torsion-free,
$M$ is a manifold. Now, $M$ is a regular cover of $O_1$,
with covering group ${\Bbb Z}/2{\Bbb Z} \times {\Bbb Z}/2{\Bbb Z}$.
Hence, each arc and circle of ${\rm sing}(O_1)$ has order 2.
Since $\pi_1(O_1)$ contains ${\Bbb Z}/2{\Bbb Z} \times {\Bbb Z}/2{\Bbb Z}$,
$O_1$ contains at least one ${\Bbb Z}/2{\Bbb Z} \times {\Bbb Z}/2{\Bbb Z}$
vertex $v$. In fact, we will show that it contains at least $4^{n-1}$
vertices. By elementary covering space theory, 
the group of covering transformations of the
cover $O_1 \rightarrow O$ equals 
$N(H)/H$, where $N(H)$ is the normaliser of $H$ in
$\prod_{i=1}^n \PSL(2,p_i)$. We claim that this group has
order at least $4^{n-1}$. To prove the claim, note that
$(I,I, \dots, I, A_i, I, \dots, I)$ and
$(I,I, \dots, I, B_i, I, \dots, I)$ both commute
with $(A_1, \dots, A_n)$ and $(B_1, \dots, B_n)$.
In particular, they lie in $N(H)$. The group these elements
generate has order $2^{2n}$. Hence, $N(H)/H$ has 
order at least $4^{n-1}$. No covering transformation
can fix $v$, because the local group of each singular point
of ${\rm sing}(O)$ does not contain ${\Bbb Z}/2{\Bbb Z} \times {\Bbb Z}/2{\Bbb Z}$
as a proper subgroup. Hence the orbit of $v$ under the
group of covering transformations has order at least $4^{n-1}$. In particular,
$O_1$ contains at least this many vertices. 
Since $n$ was an arbitrary positive integer, the theorem
is proved. $\square$

\noindent {\bf Lemma 4.7.} {\sl Let $X$ be a
finite trivalent graph with $V$ vertices.
Then $X$ contains a simple closed curve with
at most $2\log_2((V+2)/3) + 2$ edges.}

\noindent {\sl Proof.} We may assume that $X$ is connected. Give $X$ the path metric where
each edge has length 1. For any vertex $v$, let 
$R_1(v)$ be the minimal radius of a ball centred at $v$
that is not a tree. Fix a vertex $v$ where $R_1(v)$ has minimal value, 
and set $R = \lceil R_1(v) \rceil$.
Then the ball of radius $R$ around $v$ contains a simple
closed curve of length at most $2R$. We claim that $$R \leq \log_2((V+2)/3) + 1.$$

For any non-negative integer $r$, let $B(r)$ be the
ball of radius $r$ around $v$. So, $B(R - 1)$ is a tree. 
The number of vertices in this
tree is equal to $3(2^{R-1}-1)+1$. This is a lower bound
for $V$. So, 
$$3(2^{R - 1} -1)+1 \leq V$$
and therefore
$$R \leq \log_2((V+2)/3) + 1.$$
$\square$

\noindent {\sl Proof of Theorem 4.1.}
By Theorem 4.5, $O$ has a finite cover $O_1$ 
such that each arc and circle of ${\rm sing}(O_1)$
has order 2 and which
contains at least $50$ vertices. Starting with
$O_1$, we will construct a tower of finite covers
$\{ O_i \}$. Let $n_i$ be the number of vertices of $O_i$.
We will ensure that the following inequality
holds for each $i$:
$$n_{i+1} \geq 2 n_i - 4(\log_2((n_i+2)/3) + 1)\eqno{(\ast)}$$
Suppose that $H_1(|O_i|; {\Bbb Z}/2{\Bbb Z})$ is non-trivial.
Then $|O_i|$ has a 2-fold cover $|O_{i+1}|$,
with underlying orbifold $O_{i+1}$. Clearly,
the number of vertices is doubled, and so $(\ast)$ holds.
So, suppose that $|O_i|$ is a mod 2 homology 3-sphere.
Using Lemma 4.7, pick a simple closed
curve $C$ in ${\rm sing}_2(O_i)$ with length
at most $2\log_2((n_i+2)/3) + 2$. This bounds a
compact embedded surface $S$. By a small isotopy,
we may assume that $S$ intersects
${\rm sing}_2(O_i)$ in $\partial S$ and in a finite
number of points in the interior of $S$. Let $O_{i+1}$
be the 2-fold cover of $O_i$ dual to $S$.
Then each vertex in ${\rm sing}_2(O_i) - C$ has inverse
image equal to 2 vertices in $O_{i+1}$.
Thus, $(\ast)$ holds.

We claim that when $n_i$ is a sequence satisfying
$(\ast)$ and where $n_1 \geq 50$, then
$$\inf_i n_i/2^i > 0.$$
To prove this, we will establish the following inequality,
by induction on $i$:
$$n_i \geq 2^{i} \left ( 1+{24 \over i} \right ).$$
This holds for $n_1$ by our hypothesis that $n_1 \geq 50$.
To prove the inductive step, note that
$$\eqalign{
n_{i+1} &\geq 2n_i - 4(\log_2((n_i+2)/3) + 1) \cr
&\geq 2n_i - 4\log_2 n_i \cr
&\geq 2^{i+1} \left (1+{ 24 \over i} \right) - 4 \left (i + \log_2
\left (1+{ 24 \over i} \right) \right ) \cr
&\geq 2^{i+1} \left (1+{ 24 \over i} \right) - 4 (i + 5) \cr
&\geq 2^{i+1} \left (1+{ 24 \over i+1} \right).}$$
The second inequality holds because $n_i \geq 4$.
The third is true because $2x - 4 \log_2 x$ is an increasing
function of $x$ when $x > 2/\log 2$. The final
inequality holds because
$${24 \over i} - {i+5 \over 2^{i-1}} \geq {24 \over i+1} \Leftrightarrow
{24 \over i(i+1)} \geq {i +5\over 2^{i-1}},$$
which certainly holds for all integers $i \geq 1$.
So,
$${n_i \over 2^i} \geq \left( 1+ {24 \over i} \right )$$
which has positive infimum. Thus, $\{ O_i\}$ has
linear growth of mod 2 homology. $\square$
\vfill\eject

\centerline{\caps 5. Largeness criteria}
\vskip 6pt

The main result of this section is a largness criterion
for certain hyperbolic 3-orbifolds. The next
theorem is the starting point for \S 6.

\noindent {\bf Theorem 5.1.} {\sl Let $O$ be a 3-orbifold 
(with possibly empty singular locus) commensurable
with a closed orientable hyperbolic 3-orbifold that contains ${\Bbb Z}/2{\Bbb Z} \times {\Bbb Z}/2{\Bbb Z}$
in its fundamental group. Suppose that $vb_1(O) \geq 4$. Then $\pi_1(O)$
is large.}

The proof of this is independent
of the material in \S 4. The principal ingredient is the following result
(Theorem 1.2 of [16])

\noindent {\bf Theorem 5.2.} {\sl Let $G$ be a
finitely presented group, and suppose that,
for each natural number $i$, there is a
triple $H_i \geq J_i \geq K_i$ of finite index
normal subgroups of $G$ such that
\item{(i)} $H_i/J_i$ is abelian for all $i$;
\item{(ii)} $\lim_{i \rightarrow \infty} 
((\log [H_i : J_i]) / [G:H_i]) = \infty$;
\item{(iii)} $\limsup_i (d(J_i/K_i) / [G:J_i])  > 0$.

\noindent Then $K_i$ admits a surjective homomorphism
onto a free non-abelian group, for infinitely
many $i$.}

Here, $d( \ )$ is the minimal number of generators of a group.

We will need the following corollary.

\noindent {\bf Corollary 5.3.} {\sl Let $G$ be a finitely
presented group, and let $\phi \colon G \rightarrow {\Bbb Z}$
be a surjective homomorphism. Let $G_i$ be
$\phi^{-1}(i {\Bbb Z})$. Suppose that, for some
prime $p$, $\{ G_i \}$ has linear growth of
mod $p$ homology. Then $G$ is large.}

\noindent {\sl Proof.} Set $H_i = G$, set $J_i = G_i$
and let $K_i = [G_i,G_i]G_i^p$. Then it is trivial
to check that the conditions of Theorem 5.2 hold.
$\square$

A key hypothesis in Corollary 5.3 is linear growth of
mod $p$ homology. The following gives a situation where
this is guaranteed to hold.

\noindent {\bf Proposition 5.4.} {\sl Let $O$ be a compact
orientable 3-orbifold, and let $C$ be a 
component of ${\rm sing}^0_p(O)$ for some prime $p$.
Let $p_i \colon |O_i| \rightarrow |O|$ ($i\in {\Bbb N}$) be 
distinct finite covering spaces of $|O|$ such that
the restriction of $p_i$ to each component of
$p_i^{-1}(C)$ is a homeomorphism onto $C$.
Let $O_i$ be the corresponding covering spaces
of $O$. Then $\{ O_i \}$ has linear growth of mod $p$ homology.}

\noindent {\sl Proof.} By Proposition 3.1, we have
$$d_p(O_i) \geq |{\rm sing}_p^0(O_i)| \geq
[O_i :O]. \ \square$$ 

Combining Corollary 5.3 and Proposition 5.4, we have the following.

\noindent {\bf Theorem 5.5.} {\sl Let $O$ be a 
compact orientable 3-orbifold. Suppose that
$\pi_1(O)$ admits a surjective homomorphism $\phi$
onto ${\Bbb Z}$, and that some component of ${\rm sing}_p^0(O)$
has trivial image under $\phi$, for some prime $p$.
Then $\pi_1(O)$ is large.}

\noindent {\sl Proof.} Each meridian of the
singular locus of $O$ represents a torsion element
of $\pi_1(O)$. Hence its image under $\phi$ is
trivial. Thus, $\phi$ factors through a homomorphism
$\psi \colon \pi_1(|O|) \rightarrow {\Bbb Z}$.
Let $|O_i|$ be the covering space of $|O|$
corresponding to $\psi^{-1}(i {\Bbb Z})$, and
let $O_i$ be the corresponding cover of
$O$. Proposition 5.4 gives that 
$\{ O_i \}$ has linear growth of mod $p$ homology. Thus, by Corollary 5.3
(with $G = \pi_1(O)$ and $G_i = \pi_1(O_i)$),
$\pi_1(O)$ is large.
$\square$

\noindent {\bf Remark 5.6.} Suppose that the singular
locus of $O$ contains a circle component and
that $b_1(O) \geq 2$. Then such a homomorphism
$\phi$ as in Theorem 5.5 may always be
found.

\noindent {\sl Proof of Theorem 5.1.} By hypothesis, $O$ has a finite
cover $O'$ such that $b_1(O') \geq 4$. Let $O''$ be
the hyperbolic orbifold, commensurable with $O$, containing ${\Bbb Z}/2{\Bbb Z} \times {\Bbb Z}/2{\Bbb Z}$
in its fundamental group. Now,
$O'$ and $O''$ are commensurable, and hence
they have a common cover $O'''$, say. Since
$O'''$ is hyperbolic, it has a manifold cover $M$.
We may assume that $M$ regularly covers $O''$. Now,
$b_1$ does not decrease under finite covers,
and so $b_1(M) \geq 4$. Since $M \rightarrow O''$
is a regular cover, it has a group of covering
transformations $\pi_1(O'')/\pi_1(M)$.
This group acts on the manifold $M$
with quotient $O''$. Now, $\pi_1(O'')$ contains
${\Bbb Z}/2{\Bbb Z} \times {\Bbb Z}/2{\Bbb Z}$, and hence some singular
point of $O''$ has local group that contains
${\Bbb Z}/2{\Bbb Z} \times {\Bbb Z}/2{\Bbb Z}$. The group of
covering transformations must contain the local group
of this vertex. Hence,
$\pi_1(O'')/\pi_1(M)$ contains ${\Bbb Z}/2{\Bbb Z}
\times {\Bbb Z}/2{\Bbb Z}$. Let $h_1$ and $h_2$ be
the commuting covering transformations of $M$ corresponding to the
generators of ${\Bbb Z}/2{\Bbb Z} \times {\Bbb Z}/2{\Bbb Z}$.
These are involutions.
Let $h_3$ be the composition of $h_1$ and $h_2$,
which also is an involution. 
For $i=1$, 2 and 3, let $O_i$ be the quotient $M/h_i$.
Since $h_i$ has non-empty fixed point set,
${\rm sing}(O_i)$ is a non-empty collection of
simple closed curves with order 2.

We claim that, for at least one $i \in \{1, 2,  3 \}$,
$b_1(O_i) \geq 2$. 
Now, each $h_i$ induces an automorphism $h_{i \ast}$
of $H_1(M; {\Bbb R})$. Note that $h_{i \ast}$ is diagonalisable
because its minimum polynomial divides $x^2 -1$ and so splits
as a product of distinct linear factors. Thus, $H_1(M; {\Bbb R})$
decomposes as a direct sum of eigenspaces of $h_{i \ast}$.
It is clear that $b_1(O_i)$ is equal to the
dimension of the $+1$ eigenspace of $h_{i \ast}$
(see Proposition III.10.4 in [2] for example).
Suppose that this is at most 1 for
$i = 1$ and $2$. Then the dimension of
the $-1$ eigenspace is at least 3 for
$i = 1$ and $2$. Hence, the intersection of
these eigenspaces has dimension at least 2.
This lies in the $+1$ eigenspace for
$h_{3 \ast}$, and so $b_1(O_3) \geq 2$,
proving the claim.
So, by Theorem 5.5 and Remark 5.6,
$\pi_1(O_i)$ is large and hence so is $\pi_1(O)$. $\square$

\vskip 18pt
\centerline{\caps 6. Largeness for arithmetic hyperbolic 3-orbifolds}
\vskip 6pt

An easy consequence of Theorem 1.2 and Theorem 5.1 is the following.

\noindent {\bf Theorem 6.1.} {\sl Let $M$ be an arithmetic 3-manifold.
Suppose that $vb_1(M) \geq 4$. Then $\pi_1(M)$ is large.}

\noindent {\sl Proof.}
When $M$ is closed this is immediate from Theorem 1.2 and Theorem 5.1.
When $M$ has non-empty boundary, the result follows from [5],
which shows that any non-compact finite-volume hyperbolic 3-manifold
has large fundamental group. $\square$

The hypothesis $vb_1(M) \geq 4$ is known to hold in various circumstances
as we now discuss.

The first situation is:

\noindent {\bf Corollary 6.2.} {\sl Let $\Gamma$ be an arithmetic
Kleinian group, with the property that some congruence subgroup has
positive $b_1$. Then $\Gamma$ is large.}

\noindent {\sl Proof.} By Borel's theorem (Theorem 2.7) if a congruence
subgroup has positive first betti number then we have $vb_1(\Gamma) = \infty$.
Hence we can deduce largeness from Theorem 6.1. $\square$

As a particular case of this we have

\noindent {\bf Corollary 6.3.} {\sl Let $\Gamma$ be an arithmetic
Kleinian group, with the property that $vb_1(\Gamma) >0$ by arithmetic
methods. Then $\Gamma$ is large.}

We discuss some particular examples of this situation at the end of this
section.

Our next result allows us to weaken the ${\Bbb Z}/2{\Bbb Z} \times {\Bbb Z}/2{\Bbb Z}$ and 
$vb_1\geq 4$ assumption.

\noindent {\bf Theorem 6.4.} {\sl Let $M$ be an arithmetic hyperbolic 
3-manifold commensurable with an orbifold $O={\Bbb H}^3/\Gamma$ 
such that $ \Gamma$ contains $A_4$, $S_4$ or $A_5$.
Then $\pi_1(M)$ is large.}

To prove this theorem we need to recall the following result of
Clozel [4] stated in a way that is convenient for us.

\noindent {\bf Theorem 6.5.} {\sl Let $\Gamma$ be an arithmetic
Kleinian group, with invariant trace-field and quaternion algebra $k$ and $B$ respectively.
Assume that for every place $\nu\in \Ram_f(B)$, $k_\nu$ contains no quadratic extension of
${\Bbb Q}_p$ where $p$ is a rational prime and $\nu |p$. Then $\Gamma$ is commensurable with
a congruence subgroup with positive first Betti number.}

\noindent {\sl Proof of Theorem 6.4.}
We can assume that $\Gamma$ is cocompact, otherwise the result follows
from [5]. By Theorem 2.5 the 
invariant trace-field and quaternion algebra of $\Gamma$ satisfies the
conditions of Theorem 6.5. Hence we can apply Corollary 6.2 to complete
the proof. $\square$

\noindent{\bf Examples of Corollary 6.3:}

\item{1.} It is known that any arithmetic Kleinian group arising from a
quaternion algebra $B/k$ (as in \S 2.1) with $[k:k\cap {\Bbb R}]=2$
have congruence covers with $vb_1>0$ (see [14], [20] or [27]). Hence
these are large.

\item{2.} In [35] it is shown that if $k$ has one complex place and 
$[k:{\Bbb Q}]\leq 4$ then any arithmetic Kleinian group arising from
an algebra $B/k$ satisfies the hypothesis of Corollary 6.3 and hence is
large.

\item{3.} As a particular case of 2 above, let $M_W$ denote
the Weeks manifold, 
the smallest arithmetic hyperbolic 3-manifold.
>From [3], $M_W$
has invariant trace field of degree 3. Hence $\pi_1(M_W)$ is large.

\item{4.} Let $\Sigma$ denote an arithmetic integral homology 3-sphere. 
The invariant trace-field of $\Sigma$ (denoted $k$)
has even degree over ${\Bbb Q}$,
and the invariant quaternion algebra of $\Sigma$ is 
unramified at all finite places (see for example [32] Theorem 6.4.3). 
Hence Clozel's theorem (Theorem 6.5)
applies to prove $vb_1(\Sigma)>0$, and Corollary 6.2 applies to prove
largeness.  

\item{} As an example of an arithmetic integral homology 3-sphere
one can take the 3-fold cyclic branched cover of the
$(-2,3,7)$-pretzel knot. The invariant trace field is ${\Bbb
Q}(\theta)$ where $\theta$ has minimal polynomial $x^6 - x^5 - x^4 +
2x^3 - 2x^2 - x + 1$. This generates a field of signature $(4,1)$ and
discriminant $-104483$. All of this can be checked using Snap (see [6] for
a discussion of this program).

\noindent We close this section with an example of a
commensurability class of arithmetic 3-orbifolds for which no method
currently known applies to provide a cover with positive first Betti
number.

Let $p(x)=x^5-x^3-2x^2+1$. Then $p$ has three real roots and one pair
of complex conjugate roots. Let
$t$ be a complex root and let $k={\Bbb Q}(t)$. 
Now $k$ has one complex place and its Galois group is $S_5$. 
There is a unique prime
$\cal P$ of norm $11^2$ in $k$. It follows that $k_{\cal P}$ is a quadratic
extension of ${\Bbb Q}_{11}$. Take $B$ ramified at the real 
embeddings and the prime $\cal P$. Then it is unknown whether any arithmetic
Kleinian group arising from $B$ has a cover with positive first Betti number.

Briefly, if $\Gamma$ is any group in the commensurability class, then
since $k$ has odd degree, there are no non-elementary Fuchsian subgroups 
(see [32] Chapter 9). The result of Clozel (see Theorem 6.5)
does not apply by the condition on $\cal P$, and none of the papers
[14], [20] or [35] apply since $[k:{\Bbb Q}]=5$ and the 
Galois group is $A_5$.

\vskip 18pt
\centerline{\caps 7. Largeness and Property $(\tau)$}
\vskip 6pt

We begin by recalling the definition of Property $(\tau)$.

\noindent {\bf Definition.} Let $X$ be a finite graph,
and let $V(X)$ denote its vertex set. For any subset $A$ of
$V(X)$, let $\partial A$ denote those edges with one endpoint
in $A$ and one not in $A$. Define the {\sl Cheeger constant}
of $X$ to be
$$h(X) = \min \left \{ {|\partial A| \over |A|} :
A \subset V(X) \hbox{ and } 0 < |A| \leq |V(X)|/2 \right \}.$$
Let $G$ be a group with a finite generating set $S$.
For any subgroup $G_i$ of $G$, let $X(G/G_i;S)$ be the
Schreier coset graph of $G/G_i$ with respect to $S$.
Then $G$ is said to have {\sl Property $(\tau)$ with
respect to a collection of finite index subgroups} $\{ G_i \}$
if $\inf_i h(X(G/G_i; S)) > 0$. This turns out not to depend on
the choice of finite generating set $S$.  Also, $G$ is said to have 
{\sl Property 
$(\tau)$} if it has Property $(\tau)$ with respect to the collection of
all subgroups of finite index in $G$.

Lubotzky and Sarnak have made the following conjecture.

\noindent {\bf Conjecture 7.1.} (Lubotzky-Sarnak) {\sl
The fundamental group of any closed hyperbolic 3-manifold
does not have Property $(\tau)$.}

A slight variant is the following.

\noindent {\bf Conjecture 7.2.} {\sl For any closed orientable
3-manifold $M$ with infinite fundamental group, $\pi_1(M)$ does
not have Property $(\tau)$.}

It is a fairly routine argument that if we assume
the solution to the Geometrisation Conjecture,
then Conjectures 7.1 and 7.2 are equivalent. We sketch
this argument in an appendix.
%
%
%

Now it is well-known that having $vb_1>0$ implies the
Lubotzky-Sarnak conjecture. For if $M$ is a hyperbolic 3-manifold,
and we assume that $M$ is finitely covered by a 3-manifold $\tilde M$
with $b_1(\tilde M) > 0$, then $\pi_1(\tilde M)$ admits a surjective
homomorphism $\phi$ onto ${\Bbb Z}$. Let $G_i$ be $\phi^{-1}(i{\Bbb Z})$.
Then it is not hard to prove that $\pi_1(M)$ does not have Property
$(\tau)$ with respect to $\{ G_i \}$.

The main result in this section is as follows.

\noindent {\bf Theorem 7.3.} {\sl Conjecture 7.2 implies that 
any lattice in ${\rm PSL}(2, {\Bbb C})$ that contains
${\Bbb Z}/2{\Bbb Z} \times {\Bbb Z}/2{\Bbb Z}$ is
large.} 

Theorem 7.3 has the following surprising conclusion for arithmetic Kleinian
groups.

\noindent {\bf Corollary 7.4.} {\sl Conjecture 7.2  implies that
any arithmetic Kleinian group is large.}
%
%
%

There are two parts to the proof of Theorem 7.3:
\item{1.} Let $O$ be the orbifold quotient of ${\Bbb H}^3$ by the given
lattice in ${\rm PSL}(2, {\Bbb C})$.
Prove that $O$ has a finite cover $\tilde O$ where
$|\tilde O|$ has infinite fundamental group, and where
${\rm sing}_2^-(\tilde O)$ is non-empty (Theorem 8.1).
\item{2.} Prove a result (Theorem 9.1) analogous to Corollary 5.3, which
states that $\pi_1(\tilde O)$ is large, provided ${\rm sing}^-_2(\tilde O)$ is non-empty
and $\pi_1(|\tilde O|)$ does not have Property $(\tau)$.

%
%
%
%

\vskip 18pt
\centerline {\caps 8. An underlying manifold with infinite
fundamental group}
\vskip 6pt

The main theorem in this section is the following.

\noindent {\bf Theorem 8.1.} {\sl Let $O$ be a closed orientable
hyperbolic 3-orbifold that contains ${\Bbb Z}/2{\Bbb Z} \times
{\Bbb Z}/2{\Bbb Z}$ in its fundamental group. Then $O$ has a
finite cover $\tilde O$ such that ${\rm sing}_2^-(\tilde O)$ is non-empty, 
and where $|\tilde O|$ admits an infinite tower of finite covers
$\{ |O_i| \rightarrow |\tilde O| \}$. In particular,
$\pi_1(|\tilde O|)$ is infinite.}

The first step is to prove the following extension of Lemma 4.7.

\noindent {\bf Proposition 8.2.} {\sl Let $X$ be a
finite trivalent graph.
Then $X$ contains a connected subgraph $Y$ such that
$b_1(Y) = 2$ and the number of edges of $Y$ is
at most $6\log_2(b_1(X) - 1) + 12$.}

\noindent {\sl Proof.} We may assume that $X$ is connected.
Consider the path metric on $X$,
where each edge has length 1. For any vertex $v$ and non-negative integer 
$n \leq b_1(X)$, let 
$R_n(v)$ be the minimal radius of a ball centred at $v$
that contains a subgraph $Y$ with $b_1(Y) \geq n$.
Note that $b_1(X) \geq 2$ and so $R_2(v)$ is well-defined.
Fix a vertex $v$ where $R_2(v)$ has minimal value, and set $R = \lceil R_2(v) \rceil$.
We claim that $$R \leq \log_2(b_1(X) - 1) + 2.$$

There are three cases to consider: when $R_1(v) \leq R_2(v) - 1$,
when $R_1(v) = R_2(v) - 1/2$
and when $R_1(v) = R_2(v)$. Let us concentrate on the first case. 

For any non-negative integer $r$, let $B(r)$ be the
ball of radius $r$ around $v$. 
Since $R_1(v) < R_2(v)$, $b_1(B(R_1(v))) = 1$.
Thus, $B(R_1(v))$ contains a unique simple closed curve $C$
consisting of at most $2 R_1(v)$ edges. Since we are assuming
$R - 1 \geq R_1(v)$, $B(R - 1)$ is obtained from $C$ by
attaching trees. The distance between any two vertices
of $B(R - 1)$ is at most $2R - 2$. By definition of $R$,
there exist two vertices in $B(R - 1)$
joined by a path $p_1$ of length at most 2 in $X - B(R-1)$,
where $p_1$ is either an arc or simple closed curve.
Let $p_2$ be the shortest path in $B(R-1)$ from one of these
vertices to $C$. If the endpoints of $p_1$ are distinct, let $p_3$ be the shortest
path in $B(R-1)$ from the other vertex to $C \cup p_2$; otherwise
let $p_3$ be the empty set.
Set $Y = p_1 \cup p_2 \cup p_3 \cup C$.
Then it is clear that $Y$ is connected and $b_1(Y) = 2$.
By construction, $Y$ consists of at most
$6R$ edges.

Now, $C$ runs through at most two of the three
edges adjacent to $v$. Thus, if one were to
remove $v$ and its adjacent edges from $B(R - 1)$, one component would
consist of a based binary tree. The number of vertices in this
tree, together with $v$, is equal to $2^{R-1}$. This is a lower bound
for the number $V$ of vertices in $X$.
But, since $X$ is connected and trivalent, $b_1(X) = {1 \over 2}V + 1$.
So, $$2^{R - 1} \leq V = 2 b_1(X) - 2$$
and therefore
$$R \leq \log_2(b_1(X) - 1) + 2.$$
This proves the claim when $R_1(v) \leq R_2(v) -1$.
In the remaining cases, the proof is similar
but simpler, and so is omitted.

We have already seen that $Y$ contains at most
$6R$ edges, and so the proposition is proved.
$\square$

\noindent {\bf Lemma 8.3.} {\sl Let $O$ be a closed
orientable 3-orbifold with non-empty singular locus. 
Then $\pi_1(O)$ has a finite presentation
$\langle X | R \rangle$ with 
$$|R| - |X| \leq 2 b_1({\rm sing}(O)) - 2.$$}

\noindent {\sl Proof.} Let $M$ be the manifold 
obtained from $O$ by removing an open regular neighbourhood
of its singular locus. Then $\pi_1(M)$ has a presentation
$\langle X | R' \rangle$ with $|R'| - |X|  
= {1 \over 2}\chi(\partial M) - 1 = \chi({\rm sing}^-(O)) - 1$. We obtain
$\pi_1(O)$ from $\pi_1(M)$ by adding relations
that are powers of the meridians of the singular locus.
For each circle component of the singular locus,
there is one such meridian. For each graph component $Y$,
the number of meridians is equal to $-3 \chi(Y)$.
Hence, $\pi_1(O)$ has a presentation
$\langle X | R \rangle$ with 
$$\eqalign{
|R| - |X| &= -2\chi({\rm sing}^-(O)) + |{\rm sing}^0(O)| -1 \cr
&= 2b_1({\rm sing}^-(O)) - 2|{\rm sing}^-(O)| + 2b_1({\rm sing}^0(O))
- |{\rm sing}^0(O)| - 1 \cr
&\leq 2b_1({\rm sing}(O)) - 2.}$$
$\square$

The final step in the proof of Theorem 8.1 is the
following proposition.

\noindent {\bf Proposition 8.4.} {\sl Let $O$ be a 
closed orientable 3-orbifold such that every arc and circle of the singular
locus has singularity order 2.
Let $Y$ be a connected subgraph of ${\rm sing}(O)$ that has
$b_1(Y) = 2$ and at most $6\log_2(b_1({\rm sing}(O)) - 1) + 12$
edges. Let $G$ be the fundamental group of $O$, and
let $K$ be the subgroup normally generated
by the meridians that encircle $Y$. Let $\tilde O$ be the
covering space of $O$ corresponding to the subgroup $K G^2$. 
Suppose that $d_2(O) \geq 81$. Then $|\tilde O|$ has an infinite
tower of finite covers $|O_i| \rightarrow |\tilde O|$. In particular,
$\pi_1(|\tilde O|)$ is infinite.}

\noindent {\sl Proof.} 
Let $M_1$ be the meridians of $Y$, one for each edge of ${\rm sing}(O)$ that
lies in $Y$. Let $M_2$ be those meridians
of ${\rm sing}(O)$ that lie in $KG^2$
but that are not meridians of $Y$. 
Consider the group $\Gamma=\pi_1(O)/\langle \! \langle M_1, M_2 \rangle \! \rangle$.
This is equal to the fundamental group of an orbifold $O'$
with the same underlying manifold as $O$, and with singular
set that is a subgraph of ${\rm sing}(O)$.
Lemma 8.3 states that $\pi_1(O')$ has a finite presentation $\langle X | R \rangle$
where the number of relations minus the number of generators is
at most $2b_1({\rm sing}(O')) -2$, which is at most
$2b_1({\rm sing}(O)) -2$. By Proposition 3.1,
this is at most $2 d_2(O) - 2$. 
Now, adding the relations in $M_1$ to $\pi_1(O)$ reduces
$d_2$ by at most $|M_1|$. Then adding the relations in
$M_2$ does not affect $d_2$, because they lie in $KG^2$.
Thus, $d_2(\Gamma) \geq d_2(O) - |M_1|$.
Hence, 
$$d_2(O) \geq d_2(\Gamma) \geq d_2(O) - |M_1| \geq d_2(O) - 6\log_2(d_2(O) - 1) - 12.$$
Therefore,
$$\eqalign{
&d_2(\Gamma)^2/4 - |R| + |X| - d_2(\Gamma)\cr
&\qquad \geq d_2(\Gamma)^2/4 - 2 d_2(O) + 2 - d_2(\Gamma) \cr
&\qquad \geq (d_2(O) - 6\log_2(d_2(O) - 1) - 12)^2/4 - 3 d_2(O) + 2\cr
&\qquad > 0.}$$
The last inequality is a consequence of the assumption that
$d_2(O) \geq 81$; it is easy to check that the given function
of $d_2(O)$ is positive in this range.
Hence, by the Golod-Shafarevich theorem, $\Gamma$ has an infinite nested
sequence of finite index subgroups and hence is
infinite. The same is therefore true for $\Gamma^2$, because it has
finite index in $\Gamma$.

Let $\tilde O$ be the covering space of $O$ corresponding
to $K G^2$. We claim that $\pi_1(|\tilde O|) \cong \Gamma^2$.
This will prove that $\pi_1(|\tilde O|)$ has an infinite nested sequence
of finite index subgroups.
Now, $\pi_1(|\tilde O|)$ is obtained from $\pi_1(\tilde O)
= K G^2$ by quotienting each meridian of the singular
locus. But the meridians in $O$ that lift to
meridians of singular components of $\tilde O$ are
precisely those lying in $\langle \! \langle M_1, M_2 \rangle \! \rangle$. 
Hence, $\pi_1(|\tilde O|)
\cong K G^2 / \langle \! \langle M_1, M_2 \rangle \! \rangle
\cong \Gamma^2$. $\square$

\noindent {\sl Proof of Theorem 8.1.} 
By Theorem 4.5, $O$ is finitely covered by
an orbifold $O'$ such that each arc and circle of ${\rm sing}(O')$
has order 2 and where $b_1({\rm sing}_2^-(O')) \geq 81$.
By Proposition 3.1, $d_2(O') \geq 81$.
Proposition 8.2 states that ${\rm sing}_2^-(O')$ contains
a connected subgraph $Y$ such that $b_1(Y) = 2$ and which has at most $6\log_2(b_1({\rm sing}_2^-(O')) - 1) + 12$ edges,
which is at most $6\log_2(b_1({\rm sing}_2(O')) - 1) + 12$.
Let $\tilde O$ be the covering space of $O'$
corresponding to $K G^2$, where $G = \pi_1(O')$
and $K$ is the subgroup normally generated
by the meridians of $Y$. Since $\pi_1(\tilde O)$ contains $K$,
the inverse image of each edge of $Y$ in $\tilde O$ is a disjoint
union of copies of that edge. In particular,
${\rm sing}_2^-(\tilde O)$ is non-empty. By
Proposition 8.4, $|\tilde O|$ has an infinite tower of
finite covers. In particular, $|\tilde O|$ has infinite
fundamental group. $\square$

\vskip 18pt
\centerline{\caps 9. The Lubotzky-Sarnak conjecture and the largeness of orbifolds}
\vskip 6pt

In this section, we prove the following result. Together with Theorem 8.1,
this will complete the proof of Theorem 7.3.

\noindent {\bf Theorem 9.1.} {\sl Let $\tilde O$ be a compact
orientable 3-orbifold such that ${\rm sing}_2^-(\tilde O)$ is
non-empty. Let $\{ |O_i| \rightarrow |\tilde O| \}$ be a sequence of
finite-sheeted covering spaces. Suppose that $\pi_1(|\tilde O|)$ does
not have Property $(\tau)$ with respect to $\pi_1(|O_i|)$.  Then
$\pi_1(\tilde O)$ is large.}

\noindent {\sl Proof.} By passing to further finite-sheeted
covers if necessary, we may assume that each $\pi_1(|O_i|)$ is
a finite index normal subgroup of $\pi_1(|\tilde O|)$.

We are assuming that ${\rm sing}_2^-(\tilde O)$ is non-empty. Suppose
that ${\rm sing}_2^-(\tilde O)$ contains a vertex with valence 1.
(This may happen, for example, if this is a vertex of ${\rm sing}(\tilde O)$
with local group $A_5$.) If so, remove this vertex and the adjacent
edge, forming a graph $\Gamma$. Now repeat this procedure
if $\Gamma$ has a valence 1 vertex, and continue until
every vertex of $\Gamma$ has valence at least two. Remove
the components of $\Gamma$ which are circles. If a vertex has
valence 2, amalgamate its two adjacent edges into a single
edge. Repeat until every vertex of $\Gamma$ has valence at
least 3. Note that $\Gamma$ still has negative Euler characteristic.
In particular, it is non-empty. Let $\Gamma_i$ be the inverse image
of $\Gamma$ in $O_i$.

Pick a 1-vertex triangulation of $|\tilde O|$. For convenience,
we may arrange that the vertex of this triangulation is a vertex
of $\Gamma$. Place a path metric on
$|\tilde O|$ so that each edge of the triangulation has length 1 and each 3-simplex
is a regular Euclidean tetrahedron. 

The edges of this triangulation, when oriented, form a set $S$ of
generators for $\pi_1(|\tilde O|)$. We are assuming that
$\pi_1(|\tilde O|)$ does not have Property $(\tau)$ with
respect to $\{ \pi_1(|O_i|) \}$. Hence, the Cheeger constants of
the corresponding Schreier coset diagrams tend to zero.
But each such graph $X_i$ is just the 1-skeleton of $|O_i|$.
Let $A_i$ be a non-empty set of vertices in $X_i$ such that
$|\partial A_i|/|A_i| = h(X_i)$ and $|A_i| \leq |V(X_i)|/2$.
We will use $\partial A_i$ to construct a surface $S_i$ that
separates $O_i$ into two pieces $B_i$ and $C_i$. Place a 
0-cell of $S_i$ at the midpoint of each edge of $\partial A_i$.
If a 2-simplex of the triangulation of $|O_i|$ intersects $\partial A_i$, 
it does so in precisely two points. Insert into this
2-simplex a geodesic joining these two points, forming a 1-cell of $S_i$. Then,
the boundary of each 3-simplex intersects these arcs
in either the empty set, or a normal curve of length three
or four. In each 3-simplex of $|\tilde O|$ pick a representative
disc spanning the 7 different curves of length 3 and 4.
Use lifts of these discs to $|O_i|$ to construct the 2-cells
of $S_i$. This therefore defines the surface $S_i$. It divides
$O_i$ into two 3-orbifolds $B_i$ and $C_i$, say, which contain the vertices
$A_i$ and $A_i^c$ respectively.

Now, we may arrange that the singular set of $\tilde O$
is transverse to the representative normal discs in $|\tilde O|$.
Thus, ${\rm sing}(O_i)$ is transverse to $S_i$. By construction, there
is a uniform upper bound (independent of $i$) for the number
of intersection points between ${\rm sing}(O_i)$ and $S_i$
in any 3-simplex of $|O_i|$. We claim, that, viewing $S_i$
as a 2-orbifold, there is a uniform constant $K_1$ (independent of $i$) with the
following property:
$$d_2(S_i) \leq K_1 |\partial A_i|.\eqno(1)$$
This is because $d_2(S_i)$ is at most the sum of the
number of 1-cells of $S_i$ and its number of singular points.
We have already seen that the number of singular points of $S_i$
is bounded above by a constant times the number of 2-cells of $|S_i|$.
This is bounded above by the number of 1-cells of $|S_i|$.
Thus, it suffices to find a linear bound on this quantity in terms of
$|\partial A_i|$. It is at most the number
of 0-cells of $|S_i|$ times half the maximal valence of any 1-simplex of
$|\tilde O|$. But there is precisely one 0-cell of $|S_i|$ in
each edge of $\partial A_i$. This proves the claim.

We claim that there is a positive constant $K_2$, independent of $i$, such that
$$b_1(\Gamma_i \cap B_i) \geq |V(X_i)|/8 - K_2|\partial A_i|.\eqno(2)$$
Note that $\Gamma_i \cap B_i$ is a graph in which each vertex
has valence 3 or 1. The vertices with valence 1 arise at the
intersection points between $\Gamma_i$ and $S_i$.
Hence, the number of such vertices is bounded above a constant
times $|\partial A_i|$. Now, the Euler characteristic of
$\Gamma_i \cap B_i$ is equal to the sum, over all its vertices $v$,
of $1 - {\rm val}(v)/2$, where ${\rm val}(v)$ is the valence of
$v$. We arranged that the vertex of the triangulation of $|\tilde O|$
was a vertex of $\Gamma$. Hence, the number of vertices of $\Gamma_i \cap B_i$ is at
least $|A_i|$. This is at least $|V(X_i)|/4$ by Lemma 2.1 of [15].
So, for some constant $K_2$,
$$b_1(\Gamma_i \cap B_i) \geq -\chi(\Gamma_i \cap B_i) \geq |A_i|/2
- K_2|\partial A_i| \geq |V(X_i)|/8 - K_2|\partial A_i|.$$
The same inequality holds for $b_1(\Gamma_i \cap C_i)$. 

By Proposition 3.1, $d_2(B_i) \geq b_1({\rm sing}_2(B_i))
\geq b_1(\Gamma_i \cap B_i)$. Combining
this with inequalities
(1) and (2), we see that when $h(X_i)$ is small enough,
$$d_2(B_i) \geq d_2(S_i) + 2,$$
with the corresponding inequality holding also for $d_2(C_i)$.
The constant 2 here could have been replaced by any real number.
Hence, the map $H^1(B_i; {\Bbb F}_2) \rightarrow H^1(S_i; {\Bbb F}_2)$ 
induced by inclusion has 
kernel with dimension at least 2. Consider the covering space $O'$
of $O_i$ corresponding to this kernel. This has degree at least 4.
The inverse image of $C_i$ is at least four copies of $C_i$.
Now, in $C_i$, there is a properly embedded compact surface
representing a non-trivial element in the kernel of 
$H^1(C_i; {\Bbb F}_2) \rightarrow H^1(S_i; {\Bbb F}_2)$.
We may pick this surface so that it is non-separating in $C_i$,
and so that it is disjoint from $S_i$. Its inverse image
in $O'$ is a non-separating surface $F$. Let $N$ be the
number of components of $F$. Then $N \geq 4$. 

We claim there is a surjective homomorphism $\phi$ from 
$\pi_1(O')$ onto $\ast^N {\Bbb Z}/2{\Bbb Z}$, the free product
of $N$ copies of ${\Bbb Z}/2{\Bbb Z}$. The copies of ${\Bbb Z}/2{\Bbb Z}$
are indexed by the components $F_1, \dots, F_N$ of $F$. 
Let $x_i$ be the non-trivial element in the $i^{\rm th}$
copy of ${\Bbb Z}/2{\Bbb Z}$. Pick a basepoint
for $O'$ away from $F$. For each element $g$ of $\pi_1(O')$,
pick a representative loop $\ell$. Make $\ell$ transverse
to $F$ via a small homotopy (keeping the endpoints of $\ell$
fixed). As one goes round the loop $\ell$, let $F_{i_1}, \dots,
F_{i_r}$ be the components of $F$ that one meets. Define
$\phi(g)$ to be $x_{i_1} \dots x_{i_r}$. It is trivial to
check that this is invariant under a homotopy of $\ell$
relative to its endpoints. (For example, it is a consequence 
of the fact that $\phi$ is the homomorphism induced by
a collapsing map from $O'$ onto the wedge of $N$ copies
of ${\Bbb R}P^2$.) Hence, this gives a well-defined
function $\phi \colon \pi_1(O') \rightarrow \ast^N {\Bbb Z}/2{\Bbb Z}$.
It is clearly a homomorphism, since concatenation
of loops leads to concatenation of words. It is also
surjective. This is because $F$ is non-separating, and so,
for any component $F_i$ of $F$, there is a loop $\ell$,
based at the basepoint, which intersects $F_i$ once and is disjoint
from the remaining components of $F$. Hence, $\phi([\ell]) = x_i$.

Because $N \geq 4$, $\ast^N {\Bbb Z}/2{\Bbb Z}$ contains a free non-abelian
group as a finite index subgroup. The inverse image of this group
in $\pi_1(O')$ also has finite index. It surjects  this
free non-abelian group. Hence, $\pi_1(O')$ is large, as therefore
is $\pi_1(O)$. $\square$

\noindent{\bf Remark.}~If, in Theorem 9.1,
we make the extra hypotheses that the covers $|O_i|
\rightarrow | \tilde O|$ are nested, and that successive covers
$|O_{i+1}| \rightarrow |O_i|$ are regular and have degree a power of
2, then this theorem would be a consequence of the following result
(Theorem 1.1 of [18]), together with Lemma 9.3 below.

\noindent {\bf Theorem 9.2.} {\sl Let $G$ be a finitely presented
group, let $p$ be a prime and suppose that $G \geq G_1 \triangleright
G_2 \triangleright \dots$ is a nested sequence of finite index
subgroups, such that $G_{i+1}$ is normal in $G_{i}$ and has index a
power of $p$, for each $i$.  Suppose that $\{ G_i \}$ has linear
growth of mod $p$ homology. Then, at least one of the following must
hold: \item{(i)} $G$ is large; or \item{(ii)} $G$ has Property
$(\tau)$ with respect to $\{G_i \}$.

}

\noindent {\bf Lemma 9.3.} {\sl Let $\tilde O$ be a 
compact orientable 3-orbifold, such that 
${\rm sing}_2^-(\tilde O)$ is non-empty.
Let $|O_i|$ be finite covering spaces
of the manifold $|\tilde O|$ and let $O_i$ be the
corresponding covering spaces of $O$. Then
$\{ O_i \}$ has linear growth of mod $2$ homology.}

\noindent {\sl Proof.} The inverse image
of ${\rm sing}_2^-(\tilde O)$ under the map
$O_i \rightarrow \tilde O$ is ${\rm sing}_2^-(O_i)$.
Hence, $\chi({\rm sing}_2^-(O_i)) = \chi({\rm sing}_2^-(\tilde O))
[O_i : \tilde O]$. Now, 
$b_1({\rm sing}_2(O_i)) \geq |\chi({\rm sing}_2^-(O_i))|$.
Thus, Proposition 3.1 implies that
$\{ O_i \}$ has linear growth of mod $2$ homology. $\square$


We conclude this section by proving the following result; this can be deduced
from Corollary 7.4, but we give a proof below that is simpler and more
direct.

\noindent{\bf Theorem 9.4.}~{\sl If $vb_1(M) > 0$ for all
closed orientable 3-manifolds $M$ with infinite fundamental group, then
$\pi_1(N)$ is large for any arithmetic 3-manifold $N$.}

Theorem 9.4 is proved using 8.1 above together with the following.

\noindent {\bf Theorem 9.5.} {\sl Let $O$ be a closed orientable hyperbolic
3-orbifold containing
${\Bbb Z}/2{\Bbb Z} \times {\Bbb Z}/2{\Bbb Z}$ in its fundamental group. 
Suppose that $vb_1(|O|) > 0$. Then $\pi_1(O)$ is large.}

\noindent {\sl Proof.} By Proposition 4.4, $O$ is finitely covered
by an orbifold $O'$ such that ${\rm sing}_2^-(O')$ is
non-empty. We are assuming that $vb_1(|O|) > 0$.
So, there is a finite cover $|O''| \rightarrow |O|$
such that $b_1(|O''|) > 0$. Consider the covering space
$O'''$ of $O$ corresponding to $\pi_1(O') \cap \pi_1(O'')$.
We claim that $O''' \rightarrow O'$ descends to a cover
$|O'''| \rightarrow |O'|$ between underlying manifolds.
To prove this, it suffices to show that each torsion element
$g$ in $\pi_1(O')$ lies in $\pi_1(O''')$. But $g$
maps to a torsion element of $\pi_1(O)$ and this must
lie in $\pi_1(O'')$ because $|O''| \rightarrow |O|$
is a cover. Thus, $g$ lies in $\pi_1(O'') \cap \pi_1(O') = \pi_1(O''')$,
as required. Because $|O'''| \rightarrow |O'|$ is a cover,
${\rm sing}_2^-(O''')$ is non-empty.
Since $b_1(O'')$ is positive, so too is $b_1(O''')$.

Let $G = \pi_1(O''')$,
and let $\phi \colon \pi_1(O''') \rightarrow {\Bbb Z}$ be a surjective
homomorphism. Let $G_i$ be
$\phi^{-1}(i {\Bbb Z})$, and let $O_i$ be the corresponding
covering space of $O'''$. By Lemma 9.3, $\{ \pi_1(O_i) \}$
has linear growth of mod 2 homology. So, by Corollary 5.3,
$\pi_1(O''')$ is large and therefore so
is $\pi_1(O)$. $\square$

%
%
%
%

\vfill\eject
\centerline{\caps 10. The congruence subgroup property}
\vskip 6pt

In this section we show how Theorem 4.1 (proved by only the methods 
of 3-manifold topology and Kleinian groups) can be used to give a
new proof of Lubotzky's result [24].

\noindent {\bf Theorem 10.1.} {\sl No arithmetic Kleinian group
has the congruence subgroup property.}

Recall that the {\sl congruence subgroup property} is said to hold for
$\Gamma$ if any finite index subgroup of $\Gamma$ is a congruence
subgroup.  Lubotzky's original proof relied heavily on the
Golod-Shafarevich inequality and the theory of $p$-adic analytic
groups.  The aim of this section to provide a more elementary proof.
Following in the spirit of [26], the idea is to compare the number of
subgroups of $\Gamma$ of a given index with the number of congruence
subgroups with that index. Therefore, for a natural number $n$, define
$s_n(\Gamma)$ and $c_n(\Gamma)$ to be the number of subgroups of
$\Gamma$ (respectively, congruence subgroups of $\Gamma$) with index
at most $n$.

Theorem 10.1 is an immediate consequence of the following two theorems,
since they imply that $s_n(\Gamma) > c_n(\Gamma)$ for infinitely many
$n$.

\noindent {\bf Theorem 10.2.} {\sl Let $\Gamma$ be a 
lattice in ${\rm PSL}(2, {\Bbb C})$
that is commensurable with a group containing ${\Bbb Z}/2{\Bbb Z} \times {\Bbb Z}/2{\Bbb Z}$.
Then, there is a constant $k >1$ such that
$$s_n(\Gamma) \geq k^n$$
for infinitely many $n$.}

\noindent {\bf Theorem 10.3.} {\sl Let $\Gamma$ be an arithmetic Kleinian
group. Then, there is a positive constant $b$ such that
$$c_n(\Gamma) \leq n^{b \log n / \log \log n},$$
for all $n$.}

\noindent{\sl Proof of Theorem 10.2.}~Let $O$ be the orbifold
${\Bbb H}^3/\Gamma$. Suppose first that $O$ is closed.
According to Theorem 4.1, $O$ has a nested sequence of finite
covers $\{ O_i \}$ that have linear growth of mod $2$ homology.
The same conclusion holds when $O$ is non-compact, because of [4]
which guarantees 
that $\pi_1(O)$ has a finite index subgroup that has
a free non-abelian quotient. Let $\lambda$ be 
$$\inf_i d_2(O_i) / [O:O_i],$$
which is therefore positive.
Each homomorphism $\pi_1(O_i) \rightarrow {\Bbb Z}/2{\Bbb Z}$
gives a subgroup of $\pi_1(O_i)$ with index 1 or 2,
and these subgroups are all distinct. Therefore,
the number subgroups of $\pi_1(O)$ with index at most $2[O:O_i]$
is at least $2^{\lambda [O:O_i]}$. Setting $k =2^{\lambda/2}$
proves Theorem 10.2. $\square$

A proof of Theorem 10.3 for arbitrary arithmetic groups was proved by
Lubotzky in [26], and is also given
in Section 6.1 of [28]. However, given our aim of producing
a more elementary proof of Lubotzky's result in [26], we wish to avoid 
some of the technology that is used in [26] and [28].  We use these
as our guidelines but will not reproduce all of the amendments
necessary, merely commenting on salient points.

Before commencing the proof we make some comments that help
simplify some of the discussion below.  

First, following [26] and [28] we will work with the groups $\SL(2)$
rather than $\PSL(2)$. Now let $K$ be an arbitrary number field with ring of
integers $R_K$.  Then $\SL(2,R_K)$ contains the family of congruence
subgroups obtained in the usual way as:
$$\Gamma(J)=\ker(\SL(2,R_K)\rightarrow \SL(2,R_K/J)),$$
\noindent where $J\subset R_K$ is an ideal.  

Connecting with the discussion in \S 2.4, $J.M(2,R_K)$ is a proper 
2-sided integral ideal of
$M(2,R_K)$ contained in $M(2,R_K)$, and the elements $\alpha\in
\SL(2,R_K)$ such that $\alpha-1 \in J.M(2,R_K)$ correspond precisely
to the group $\Gamma(J)$.  Also, given any other maximal order ${\cal
L}\in M(2,K)$ (not $\GL(2,K)$-conjugate to $M(2,R_K)$), then ${\cal
L}^1$ can be conjugated to be a subgroup of a group $\SL(2,R_H)$
for some number field $H$ (the Hilbert Class field of $K$ will
work, see the proof of [32] Lemma 5.2.4). 
Via this, the congruence subgroups $\Gamma({\cal L}(I))$ of \S 2.4 can be
described in $\SL(2,R_H)$ using the more traditional definition given above
(see also the discussion in [32] Chapter 6.6)

For congruence subgroups of arithmetic Kleinian
groups, using an embedding of $B$ into $\SL(2,K)$ where $K$ is a splitting
field of $B$, we can use the above discussion to describe these
congruence subgroups as subgroups of congruence subgroups of $\SL(2,R_K)$ for
certain number fields $K$.

\noindent{\bf Proof of Theorem 10.3.}

Following the above discussion, it suffices to prove our result in the
following context. 

\noindent{\bf Proposition 10.4.}~{\sl Let $K$ be a number field with ring 
of integers $R_K$ and degree $d$.
Then, there is a positive constant $b$ such that
$$c_n(\SL(2,R_K)) \leq n^{b \log n / \log \log n},$$
for all $n$.}

The key result is the following ``level versus index'' result.

\noindent {\bf Proposition 10.5.} {\sl There is some constant
$c$ with the following property. For each congruence subgroup
$H$ of $\SL(2,R_K)$, $H \geq \Gamma(J)$, for some ideal $J \subset R_K$
with $N(J) \leq c [\SL(2,R_K):H]$.}

We defer comment on the proof of this theorem and complete the proof of
Theorem 10.3.
A consequence of Proposition 10.5 is that $c_n(\SL(2,R_K))$ is
at most the sum, $\sum s_n(\SL(2,R_K/J))$, where the sum is 
over all ideals $J\subset R_K$ with $N(J) \leq cn$.
This in turn is less than the sum 
$$\sum_{m=1}^{cn}s_n(\SL(2,R_K/mR_K)).$$
Thus, the question now reduces
to a count of subgroups in the finite groups $\SL(2,R_K/mR_K)$. 
We now discuss this, 
following [26] and [28].

Define the {\sl rank} of a finite group $G$ to be 
$$\rank(G) = \sup\{d(H):H \leq G\}.$$  
An easy argument (see [28] Lemma 1.2.2)
shows the total number of subgroups of a finite group $G$ is
$\leq |G|^{\rank(G)}$.  

The groups $\SL(2,R_K/mR_K)$
decompose as $\prod \SL(2,R_K/{\cal P}_j^{a_j})$, where 
$mR_K={\cal P}_1^{a_1}\ldots {\cal P}_t^{a_t}$ is a factorization into
distinct prime (ideal) powers of
the principal ideal $mR_K$. From this, the orders of these groups
can be computed [37]. Also note that there are at most $d$ (which
recall is the degree of $K$) primes in $K$
lying above a given rational prime $p$. Given this, 
an estimate of the order that suffices is $m^{3d}$.

To compute the rank we argue as follows. 
If $\cal P$ is a $K$-prime dividing the rational prime $p$ then
$R_K/{\cal P}$ is a finite extension of ${\Bbb F}_p$ the field of $p$ elements
of degree at most $d$. Notice that if $R_{\cal P}$ denotes the ${\cal P}$-adic
integers in $K_{\cal P}$ with uniformizer $\pi_{\cal P}$ 
then 
$$\SL(2,R_K/{\cal P}^a) \cong \SL(2,R_{\cal P}/\pi^aR_{\cal P})$$
and the latter are all homomorphic images of $\SL(2,R_{\cal P})$. From
[7], standard aspects of uniform pro-$p$ groups imply that these
have rank (as pro-$p$ groups) at most $3d$ (the field $K_{\cal P}$ has degree
at most $d$ over ${\Bbb Q}_p$).  Hence we deduce:

\noindent {\bf 1.}~Let $\cal P$ be a $K$-prime.
Then $\rank(\SL(2,R_K/{\cal P}^a)) \leq 3d$.

To pass to the rank of $\SL(2,R_K/mR_K)$, the argument is as in
[26] and [28] with the extra care that we are working with a number field.

Thus if $m=p_1^{b_1}\ldots p_l^{b_l}$ then as in [26], a simple
application of the Prime Number Theorem gives $l\leq \log m / \log
\log m$.  As noted above, each prime $p_i$ splits into at most $d$
$K$-primes, and so arguing as in [26] we deduce using 1:

\noindent {\bf 2.}~$\rank(\SL(2,R_K/mR_K)) \leq 3d^2 \log m / \log \log m$.

Hence using 2 and the estimate above, it follows
that the total number of subgroups of $\SL(2,R_K/mR_K)$
is at most $m^{C\log m / \log \log m}$ where $C$ depends on $K$.
The count of congruence subgroups now finishes off as in [26].

We now discuss a proof of Proposition 10.5 in our context. 
Again, as in [26] and
[28] the key assertion concerns essential subgroups of $\SL(2,R_K/J)$.
The definition here is amended from that in [26] and [28]
to work in the number field $K$.

Following [26] and [28] we say
$H < \SL(2,R_K/J)$ is called {\sl essential} if $H$ does
not contain:
$$M(I) =  \ker(\SL(2,R_K/J)\rightarrow \SL(2,R_K/I))$$
for any $I|J$ (as ideals) with $I\neq J$.

\noindent{\bf Claim.}~{\sl There exists a constant $C'>0$
(depending on $K$) such that
for every ideal $J$, every essential subgroup $H$ of $\SL(2,R_K/J)$
satisfies 
$$[\SL(2,R_K/J):H] \geq C'N(J).$$}

\noindent{\sl Proof of Claim.} We follow the argument in [26] and [28]
adapted to our setting.  If $J$ is a prime ideal of $R_K$ then 
$R_K/J$ is a field of order $q=p^t$ for some $t\leq d$.  
It is a classical result dating back to Galois (see [37] Chapter 6)
that the minimal index of a proper subgroup of $\SL(2,R_K/J)$ 
is at least $q+1$ apart from a finite number of values of $q$. 

If $J={\cal P}^a$ is the power of a prime ideal $\cal P$, then the
argument in [26] applies in exactly the same way. The arguments in
[26] appeal to Strong Approximation, but the argument for $\SL(2)$
can be handled directly (for example by the methods of
[23]). Similarly, if $J = {\cal P}_1\ldots {\cal P}_t$ is a product of
primes whose norms are powers of distinct rational primes, then the
argument of [28] applies.

Now consider the case of $J=pR_K$ where $p$ is a rational prime that splits
completely in $K$. In this case $pR_K={\cal P}_1\ldots {\cal P}_d$.
Let $M=\SL(2,R_K/{\cal P}_1)\times \ldots \times \SL(2,R_K/{\cal P}_d)$ and
assume that $H$ is an essential subgroup of $M$.  If $\pi_i$ denotes projection
on the $i$-th factor, then we claim $\pi_i(H)$ is an essential subgroup
of $\SL(2,R_K/{\cal P}_i)$; ie $\pi_i(H)$ is a proper subgroup of 
$\SL(2,R_K/{\cal P}_i)$.

Suppose not, and assume $i=1$ for convenience. Let $M_1=\SL(2,R_K/{\cal P}_1)$
and $M'=\SL(2,R_K/{\cal P}_2)\times \ldots \times \SL(2,R_K/{\cal P}_d)$, so
$M=M_1\times M'$.
Then $\pi_1(H) = M_1$ implies by properties of the direct product that
$H \cap M_1$ is normal in $M_1$.  For all but a finite number of finite fields
$\Bbb F$, $\SL(2,{\Bbb F})$ is a central extension of a finite simple
group, and so $H\cap M_1=M_1$ or is central. Both of these can be ruled
out as in [28], using the essentialness of $H$.
Hence we deduce in this case that 
$$[M:H] > \prod[\SL(2,{\cal P}_i):\pi_i(H)] > p^d,$$
which proves the claim in this case.

Now consider $J = {\cal P}_1^{e_1}\ldots 
{\cal P}_t^{e_t}$ with some $e_i>1$. Let $J' = {\cal P}_1\ldots {\cal P}_t$,
so that $N(J')<N(J)$. Let $H$ be an essential subgroup of $\SL(2,R_K/J)$,
and $H'$ denote the projection of $H$ to $\SL(2,R_K/J')$. As in [28]
$H'$ is an essential subgroup of $\SL(2,R_K/J')$ and so the index (by above)
is at least $N(J')$. Following [28] it suffices to prove that
$H\cap M(J')$ has index at least $N(J)/N(J')$ in $M(J')$ (which is a product of
$p$-groups for various primes $p$) and this 
is completed as in the last few paragraphs
of [28] pp 116--117. $\square$

\vskip 18pt
\centerline{\caps References}
\vskip 6pt

\item{1.} {\caps A. Borel,} {\sl Cohomologie de sous-groupes discrets et
repr\'esentations de}\hfill\break
{\sl groupes semi-simples}, Ast\'erisque 32-33 (1976) 73-112.

\item{2.} {\caps  K. Brown,} {\sl Cohomology of Groups}, Graduate Texts in Math.
87 (1982), Springer-Verlag.

\item{3.} {\caps T. Chinburg, E. Friedman, K. N. Jones and A. W. Reid,} {\sl
The smallest volume arithmetic hyperbolic 3-manifold},
Annali della Scuola Normale Superiore di Pisa 30 (2001) 1--40.

\item{4.} {\caps L. Clozel}, {\sl On the cuspidal cohomology of 
arithmetic subgroups
of $\SL(2n)$ and the first betti number of arithmetic 3-manifolds},
Duke Math. J. 55 (1987) 475-486.

\item{5.} {\caps D. Cooper, D. D. Long and A. W. Reid,} {\sl Essential 
closed surfaces in bounded 3-manifolds}, Journal
of The American Math. Soc. 10 (1997) 553--563.

\item{6.} {\caps D. Coulsen, O. A. Goodman, C. D. Hodgson and W. D. Neumann,}
{\sl Computing arithmetic invariants of 3-manifolds}, Experimental
J. Math. 9 (2000) 127--152

\item{7.} {\caps J. Dixon, M. du Sautoy, A. Mann, D. Segal,}
{\sl Analytic pro-$p$ groups,} 
Cambridge Studies in Advanced Mathematics, 61 (1999),
Cambridge University Press.

\item{8.} {\caps N. Dunfield, W. Thurston,}
{\sl The virtual Haken conjecture: Experiments and examples,}
Geom. Topol. 7 (2003) 399--441.

\item{9.} {\caps P. Hall}, {\sl The Eulerian functions of a group}, Quart J.
Math. 7 (1936) 134--151.

\item{10.} {\caps J. Hempel}, {\sl Residual finiteness for $3$-manifolds},
Combinatorial group theory and topology (Alta, Utah, 1984)  379--396,
Ann. of Math. Stud., 111, Princeton Univ. Press, 1987. 

\item{11.} {\caps K. N. Jones and A. W. Reid}, {\sl  Geodesic intersections in 
arithmetic hyperbolic $3$-manifolds}, Duke Math. J.  89 (1997) 75--86.

\item{12.} {\caps T. J\o rgenson}, {\sl Closed geodesics on Riemann surfaces}, 
Proc. Amer. Math. Soc. 72  (1978) 140--142.

\item{13.} {\caps S. Kojima}, {\sl Finite covers of 3-manifolds
containing essential surface of Euler characteristic $=0$,}
Proc. Amer. Math. Soc. 101 (1987) 743--747.

\item{14.} {\caps J.-P Labesse and J. Schwermer,} {\sl On liftings and cusp
cohomology of arithmetic groups}, Invent. Math. 83 (1986) 383-401.

\item{15.} {\caps M. Lackenby}, {\sl Heegaard splittings,
the virtually Haken conjecture and Property ($\tau$)},
Invent. Math. 164 (2006) 317--359.

\item{16.} {\caps M. Lackenby,} {\sl A characterisation of
large finitely presented groups}, 
J. Algebra 287 (2005) 458--473.

\item{17.} {\caps M. Lackenby,} {\sl Covering spaces of 3-orbifolds},
Duke Math J. 136 (2007) 181-203.

\item{18.} {\caps M. Lackenby,} {\sl Large groups, Property $(\tau)$
and the homology growth of subgroups}, Preprint.

\item{19.} {\caps M. Lackenby}, {\sl Some 3-manifolds and 3-orbifolds
with large fundamental group}, Proc. Amer. Math. Soc. 135 (2007) 3393-3402.

\item{20.} {\caps J. S. Li and J. J. Millson,} {\sl On the first betti 
number of a hyperbolic manifold with an arithmetic fundamental group}, Duke
Math. J. 71 (1993) 365-401.

\item{21.} {\caps D. D. Long,} {\sl Immersions and embeddings of
totally geodesic surfaces}, Bull. London Math. Soc. 19
(1987) 481--484.

\item{22.} {\caps D. D. Long, G. Niblo,}
{\sl Subgroup separability and $3$-manifold groups.}
Math. Z. 207 (1991) 209--215.

\item{23.} {\caps D. D. Long, A. W. Reid,}
{\sl Simple quotients of hyperbolic $3$-manifold groups,} 
Proc. Amer. Math. Soc. 126 (1998) 877--880.

\item{24.} {\caps A. Lubotzky}, {\sl Group presentations, $p$-adic analytic 
groups and lattices in $\SL(2,{\Bbb C})$,} Ann. Math. 118 (1983) 115--130.

\item{25.} {\caps A. Lubotzky}, {\sl Free quotients and the first Betti number 
of some hyperbolic manifolds,} Transform. Groups  1  (1996) 71--82. 

\item{26.} {\caps A. Lubotzky}, {\sl Subgroup growth and congruence subgroups},
Invent. Math. 119 (1995) 267--295.

\item{27.} {\caps A. Lubotzky}, {\sl Eigenvalues of the 
Laplacian, the first Betti number and the congruence subgroup
problem,} Ann. Math. {\bf 144} (1996) 441--452.

\item{28.} {\caps A. Lubotzky, D. Segal},
{\sl Subgroup Growth}. Progress in Mathematics, 212. 
Birkh\"auser Verlag (2003)

\item{29.} {\caps A. Lubotzky, R. Zimmer}, {\sl
Variants of Kazhdan's property for subgroups of semisimple groups}, 
Israel J. Math. {\bf 66} (1989) 289--299.

\item{30.} {\caps J. Luecke}, {\sl Finite covers of 3-manifolds containing
essential tori,} Trans. Amer. Math. Soc. 310 (1988) 381--391.

\item{31.} {\caps C. Maclachlan, G. J. Martin}, {\sl 2-generator arithmetic
Kleinian groups}, J. f\"ur die Reine und Angew. Math. {\bf 511} (1999) 95--117.

\item{32.} {\caps C. Maclachlan, A. W. Reid}, {\sl The Arithmetic 
of Hyperbolic 3-Manifolds}, Graduate Texts in Mathematics, 219, 
Springer-Verlag (2003)

\item{33.} {\caps G. Margulis}, {\sl Discrete Subgroups of Semi-simple 
Lie Groups},
Ergebnisse der Mathematik und ihr Grenzgebeite,  Springer (1991).

\item{34.} {\caps W. Narkiewicz}, {\sl Algebraic numbers}, 
Polish Scientific Publishers (1974)

\item{35.} {\caps C. S. Rajan,} {\sl On the non-vanishing of the first 
Betti number of hyperbolic three manifolds},  Math. Ann.  330  (2004) 323--329.

\item{36.} {\caps A. W. Reid}, {\sl Isospectrality and commensurability of 
arithmetic hyperbolic $2$- and $3$-manifolds},  Duke Math. J. 65 (1992),
215--228.

\item{37.} {\caps M. Suzuki,} {\sl Group Theory I}, Grundlehren der math. 
Wissen. 247, Springer-Verlag, 1980.

%
%
%
%
%

\vfill\eject
\centerline{\caps Appendix. The equivalence of conjectures 7.1 and 7.2}
\vskip 6pt

The aim of this section is to establish the relationship between the
following conjectures.

\noindent {\bf Conjecture 7.1.} (Lubotzky-Sarnak) {\sl
The fundamental group of any closed hyperbolic 3-manifold
does not have Property $(\tau)$.}

\noindent {\bf Conjecture 7.2.} {\sl For any closed orientable
3-manifold $M$ with infinite fundamental group, $\pi_1(M)$ does
not have Property $(\tau)$.}

Our aim is to prove the following.

\noindent {\bf Theorem A.1.} {\sl If the geometrisation conjecture is
true, then Conjectures 7.1 and 7.2 are equivalent.}

\noindent {\sl Proof.} Conjecture 7.2 clearly implies Conjecture 7.1
since every closed hyperbolic 3-manifold has infinite fundamental group.
The aim is to prove the converse. Therefore, let $M$ be a closed
orientable 3-manifold with infinite fundamental group.
Our aim is to show that either $vb_1(M) > 0$ or $M$ is hyperbolic.
Conjecture 7.1 then implies that $\pi_1(M)$ does not have Property $(\tau)$.
Assuming the geometrisation conjecture,
$M$ admits a decomposition into geometric pieces.

Suppose first that $M$ is a connected sum
$M_1 \sharp M_2$, say. By geometrisation, $\pi_1(M_i)$ is
residually finite and non-trivial (see [10]). In particular, it
admits a surjective homomorphism $\phi_i$ onto a non-trivial finite
group $F_i$. Since $\pi_1(M)$ is isomorphic to
$\pi_1(M_1) \ast \pi_1(M_2)$, we therefore obtain
a surjective homomorphism from $\pi_1(M)$ onto $F_1 \ast F_2$.
But, $F_1 \ast F_2$ has a free group as a finite index
subgroup. In particular, $vb_1(F_1 \ast F_2) > 0$.
Hence, $vb_1(M) > 0$.

So, we may assume that $M$ is prime. Suppose that it
contains an essential torus. Then it is known in this case 
that either $M$ is finitely covered by a torus bundle over the circle 
or $\pi_1(M)$ is large. In
particular, $vb_1(M) > 0$ (see [13], [30] and [22]).

Consider now the case where $M$ is prime and
atoroidal. By geometrisation, it is therefore either
a Seifert fibre space or hyperbolic. In the latter case,
the proof is complete. In the former case, the proof
divides according to whether the base orbifold has
positive, zero or negative Euler characteristic.
When it is positive, the manifold is covered either by
$S^3$ or $S^2 \times S^1$. The former case is impossible,
since $\pi_1(M)$ is infinite. In the latter
case, $vb_1(M) = 1 > 0$. When the base orbifold has
zero Euler characteristic, $M$ is finitely covered
by a torus-bundle over the circle. In particular,
$vb_1(M) > 0$. When the base orbifold $O$ has negative
Euler characteristic, it is hyperbolic. Hence, in this case,
$\pi_1(O)$ is large. But the Seifert fibration induces a surjective
homomorphism from $\pi_1(M)$ onto $\pi_1(O)$, and so
$\pi_1(M)$ is also large. $\square$

\vskip 18pt

\noindent Mathematical Institute, University of Oxford,\hfill\break
\noindent 24-29 St Giles', Oxford OX1 3LB, UK.\hfill\break
\noindent Department of Mathematics, U.C.S.B, CA 93106, USA.\hfill\break
\noindent Department of Mathematics, University of Texas, Austin, TX 78712, USA.

\end

%
%
%
%
%
%

\item{2.} Use Conjecture 7.2 to give a sequence of finite covers
$|O_i|$ of $|\tilde O|$ which do not have Property $(\tau)$.
Let $O_i$ be the induced orbifold cover of $\tilde O$.
The fact that ${\rm sing}_2^-(\tilde O)$ is non-empty will
imply that $\{ O_i \}$ has linear growth of mod 2 homology.

%
%
%
%
%
%
In this section we provide the relevant extensions of work in [La1] and [La2] that will be needed
to apply in the arithmetic setting.  In particular we aim to prove Theorem 1.3.  This will follow from
Theorem 3.1 which we prove below.

\noindent {\bf Theorem 3.1.} {\sl Let $O$ be a 3-orbifold 
(with possibly empty singular locus) commensurable
with a closed hyperbolic orbifold that contains ${\Bbb Z}/2{\Bbb Z} \times {\Bbb Z}/2{\Bbb Z}$
in its fundamental group. Suppose that $vb_1(O) \geq 4$. Then $\pi_1(O)$
is large.}

\noindent Deferring the proof of Theorem 3.1 for now, the proof of Theorem 1.3 is immediate from Theorem 3.1 and
Theorem 1.2. In addition, Corollary 1.4 follows from Theorem 3.1, Theorem 1.2 together with Theorem 2.6 which gaurantees $vb_1(M)\geq 4$.

To prove Theorem 1.5 we recall the following result of Clozel [C] stated in a way that is convenient for us.

\noindent {\bf Theorem 3.10.} {\sl Let $\Gamma$ be an arithmetic
Kleinian group, with invariant trace-field and quaternion algebra $k$ and $B$ respectively.
Assume that for every place $\nu\in \Ram_f(B)$ $k_\nu$ contains no quadratic extension of
${\Bbb Q}_p$ where $p$ is a rational prime and $\nu |p$. Then $\Gamma$ is commensurable with
a congruence subgroup with positive first Betti number.}

This together with Theorem 3.1, Propositions 2.5 and 2.6 prove Theorem 1.5.$\square$
\vskip 6.0pt

\noindent {\bf 3.1}~We begin with a definition.

\noindent {\bf Definition.} Let $O$ be a compact
orientable 3-orbifold.
Let ${\rm sing}(O)$ be its singular locus,
and let $|O|$ denote its underlying 3-manifold.
Let ${\rm sing}^0(O)$ and ${\rm sing}^-(O)$
denote the components of the
singular locus with, respectively, 
zero and negative Euler characteristic.
For any prime $p$,
let ${\rm sing}_p(O)$ denote the union of the
arcs and circles in ${\rm sing}(O)$
with singularity order that is a multiple of $p$.
Let ${\rm sing}_p^0(O)$ and ${\rm sing}_p^-(O)$
denote those components of ${\rm sing}_p(O)$
with zero and negative Euler characteristic.

Note that, although ${\rm sing}(O) = {\rm sing}^0(O)
\cup {\rm sing}^-(O)$, it need not be the case
that ${\rm sing}_p(O) = {\rm sing}_p^0(O)
\cup {\rm sing}_p^-(O)$

\noindent {\bf Terminology.} If $p$ is a prime, let
${\Bbb F}_p$ denote the field of order $p$.
If $X$ is a group, space or orbifold, let $d_p(X)$
be the dimension of $H_1(X; {\Bbb F}_p)$.

\noindent {\bf Proposition 3.2.} {\sl Let $O$ be a
compact orientable 3-orbifold, and let $p$ be
a prime. Then $d_p(O) \geq b_1({\rm sing}_p(O))$.}

\noindent {\sl Proof.} Let $M$ denote the
3-manifold obtained from $O$ by a removing
an open regular neighbourhood of its singular
locus. Let $\{ \mu_1, \dots, \mu_r \}$ be a collection
of meridian curves, one encircling each
arc or circle of the singular locus.
Let $n_i$ be the singularity order of the arc
or circle that $\mu_i$ encircles. Then 
$$\pi_1(O) = \pi_1(M)/ \langle \! \langle
\mu_1^{n_1}, \dots, \mu_r^{n_r} \rangle \! \rangle.$$
Hence, 
$$H_1(O; {\Bbb F}_p) = H_1(M; {\Bbb F}_p) 
/ \langle \! \langle
\mu_1^{n_1}, \dots, \mu_r^{n_r} \rangle \! \rangle.$$
Now, when $n_i$ is coprime to $p$, quotienting
$H_1(M; {\Bbb F}_p)$
by $\mu_i^{n_i}$ is the same as quotienting by
$\mu_i$. And when $n_i$ is a multiple of
$p$, then quotienting $H_1(M; {\Bbb F}_p)$ by $\mu_i^{n_i}$ 
has no effect. Thus, if we let $M'$ be the
3-manifold obtained from $|O|$ by removing
an open regular neighbourhood of ${\rm sing}_p(O)$,
then $d_p(O) = d_p(M')$.
Now, it is a well known consequence of Poincar\'e
duality that, for the compact orientable 3-manifold $M'$,
$d_p(M') \geq {1 \over 2} d_p(\partial M') = 
b_1({\rm sing}_p(O))$, as required. $\square$

\noindent {\bf Definition.} Let $X$ be a group, space or orbifold and let $p$ be
a prime.
Then a collection $\{ X_i \}$ of finite index subgroups
or finite-sheeted covers with index or degree $[X:X_i]$
is said to have {\sl linear growth of $p$-homology}
if
$$\inf_i d_p(X_i)/[X:X_i] > 0.$$

\noindent {\bf Corollary 3.3.} {\sl Let $O$ be a 
compact orientable 3-orbifold, such that 
${\rm sing}_2^-(O)$ is non-empty.
Let $|O_i|$ be the finite covering spaces
of the manifold $|O|$ and let $O_i$ be the
corresponding covering spaces of $O$. Then
$\{ O_i \}$ has linear growth of $2$-homology.}

\noindent {\sl Proof.} The inverse image
of ${\rm sing}_2^-(O)$ under the map
$O_i \rightarrow O$ is ${\rm sing}_2^-(O_i)$.
Hence, $\chi({\rm sing}_2^-(O_i)) = \chi({\rm sing}_2^-(O))
[O_i : O]$. Now, 
$b_1({\rm sing}_2(O_i)) \geq |\chi({\rm sing}_2^-(O_i))|$.
Thus, Proposition 3.2 implies that
$\{ O_i \}$ has linear growth of $2$-homology. $\square$

\noindent {\bf Corollary 3.4.} {\sl Let $O$ be a compact
orientable 3-orbifold, and let $C$ be a 
component of ${\rm sing}^0_p(O)$ for some prime $p$.
Let $p_i \colon |O_i| \rightarrow |O|$ be the 
finite covering spaces of $|O|$ such that
the restriction of $p_i$ to each component of
$p_i^{-1}(C)$ is a homeomorphism onto $C$.
Let $O_i$ be the corresponding covering spaces
of $O$. Then $\{ O_i \}$ has linear growth of $p$-homology.}

\noindent {\sl Proof.} Again, we have
$b_1({\rm sing}_p(O_i)) \geq |{\rm sing}_p^0(O_i)| \geq [O_i :O]$. $\square$

The following proposition gives a sufficient condition for
the hypotheses of Corollary 3.3 to hold in some finite
cover.

\noindent {\bf Proposition 3.5.} {\sl Let $O$ be a closed
orientable hyperbolic 3-orbifold that contains ${\Bbb Z}/2{\Bbb Z} \times {\Bbb Z}/2{\Bbb Z}$
in its fundamental group. Then $O$ is finitely covered
by a 3-orbifold $\tilde O$ such that ${\rm sing}_2^-(\tilde O)$
is non-empty.}

\noindent {\sl Proof.} Since $O$ is hyperbolic, Selberg's lemma
implies that it has a finite-sheeted regular cover that is a manifold
$M$. Let $\tilde O$ be the cover of $O$ corresponding to the subgroup
$\pi_1(M) ( {\Bbb Z}/2{\Bbb Z} \times {\Bbb Z}/2{\Bbb Z})$ of $\pi_1(O)$.
Then $M$ regularly covers $\tilde O$ with covering group
$\pi_1(\tilde O)/\pi_1(M) = ( {\Bbb Z}/2{\Bbb Z} \times {\Bbb Z}/2{\Bbb Z}) / 
(\pi_1(M) \cap ( {\Bbb Z}/2{\Bbb Z} \times {\Bbb Z}/2{\Bbb Z})) = {\Bbb Z}/2{\Bbb Z} \times {\Bbb Z}/2{\Bbb Z}$. 
Thus, any arc or circle of the singular locus of $\tilde O$ has order 2.
Hence, ${\rm sing}_2(\tilde O) = {\rm sing}(\tilde O)$.
Since $O$ is closed, ${\rm sing}(O)$ consists of
simple closed curves and trivalent graphs.
Now, $\pi_1(\tilde O)$ contains ${\Bbb Z}/2{\Bbb Z} \times {\Bbb Z}/2{\Bbb Z}$
and therefore $\tilde O$ contains a ${\Bbb Z}/2{\Bbb Z} \times {\Bbb Z}/2{\Bbb Z}$
singular vertex. Hence, some component of ${\rm sing}_2(\tilde O)$ therefore
has negative Euler characteristic. $\square$

In [La], the following theorem was proved.

\noindent {\bf Theorem 3.6.} {\sl Let $G$ be a
finitely presented group, and suppose that,
for each natural number $i$, there is a
triple $H_i \geq J_i \geq K_i$ of finite index
normal subgroups of $G$ such that
\item{(i)} $H_i/J_i$ is abelian for all $i$;
\item{(ii)} $\lim_{i \rightarrow \infty} 
((\log [H_i : J_i]) / [G:H_i]) = \infty$;
\item{(iii)} $\limsup_i (d(J_i/K_i) / [G:J_i])  > 0$.

\noindent Then $K_i$ admits a surjective homomorphism
onto a free non-abelian group, for infinitely
many $i$.}

The following corollary of Theorem 3.6 will be useful for us.

\noindent {\bf Corollary 3.7.} {\sl Let $G$ be a finitely
presented group, and let $\phi \colon G \rightarrow {\Bbb Z}$
be a surjective homomorphism. Let $G_i$ be
$\phi^{-1}(i {\Bbb Z})$. Suppose that, for some
prime $p$, $\{ G_i \}$ has linear growth of
$p$-homology. Then $G$ is large.}

\noindent {\sl Proof.} Set $H_i = G$, set $J_i = G_i$
and let $K_i = [G_i,G_i]G_i^p$. Then it is trivial
to check that the conditions of Theorem 3.6 hold.
$\square$

The following result is the main tool we use for
establishing that certain 3-orbifolds have large fundamental
group.

\noindent {\bf Theorem 3.8.} {\sl Let $O$ be a 
compact orientable 3-orbifold. Suppose that
$\pi_1(O)$ admits a surjective homomorphism $\phi$
onto ${\Bbb Z}$, and that one of
the following holds:
\item{1.} some component of ${\rm sing}_p^0(O)$
has trivial image under $\phi$, for some prime $p$; or
\item{2.} $O$ is closed and hyperbolic, and $\pi_1(O)$ contains ${\Bbb Z}/2{\Bbb Z}
\times 
{\Bbb Z}/2{\Bbb Z}$.

\noindent Then $\pi_1(O)$ is large.}

\noindent {\sl Proof.} Each meridian of the
singular locus of $O$ represents a torsion element
of $\pi_1(O)$. Hence its image under $\phi$ is
trivial. Thus, $\phi$ factors through a homomorphism
$\psi \colon \pi_1(|O|) \rightarrow {\Bbb Z}$.
Let $|O_i|$ be the covering space of $|O|$
corresponding to $\psi^{-1}(i {\Bbb Z})$, and
let $O_i$ be the corresponding cover of
$O$. In case (1), Corollary 3.4 gives that 
$\{ O_i \}$ has linear growth of $p$-homology. Thus, by Corollary 3.7,
$\pi_1(O)$ is large.

In case (2), Proposition 3.5 implies that
$O$ is finitely covered by an orbifold $\tilde O$
such that ${\rm sing}_2^-(\tilde O)$ is non-empty.
The composition of $\pi_1(\tilde O) \rightarrow 
\pi_1(O)$ and $\phi$ gives a surjective homomorphism
from $\pi_1(\tilde O)$ onto a finite index subgroup of ${\Bbb Z}$,
and hence a surjective homomorphism $\psi$ 
from $\pi_1(|\tilde O|)$ onto ${\Bbb Z}$.
Let $|O_i|$ be the covering space of $|\tilde O|$
corresponding to $\psi^{-1}(i {\Bbb Z})$, and
let $O_i$ be the corresponding cover of $\tilde O$.
Corollary 3.3 implies
that $\{ O_i \}$ has linear growth of $2$-homology.
Corollary 3.7 then gives that $\pi_1(\tilde O)$ is
large, and hence so is $\pi_1(O)$. 
$\square$

\noindent {\bf Remark} Suppose that the singular
locus of $O$ contains a circle component and
that $b_1(O) \geq 2$. Then such a homomorphism
$\phi$ as in (1) of Theorem 1.9 may always be
found.

\noindent {\bf Proof of Theorem 3.1.} By hypothesis, $O$ has a finite
cover $O'$ such that $b_1(O') \geq 4$. Let $O''$ be
the hyperbolic orbifold containing ${\Bbb Z}/2{\Bbb Z} \times {\Bbb Z}/2{\Bbb Z}$
in its fundamental group. Now,
$O'$ and $O''$ are commensurable, and hence
they have a common cover $O'''$, say. Since
$O'''$ is hyperbolic, it has a manifold cover $M$.
We may assume that $M$ regularly covers $O''$. Now,
$b_1$ does not decrease under finite covers,
and so $b_1(M) \geq 4$. Since $M \rightarrow O$
is a regular cover, it has a group of covering
transformations $\pi_1(O'')/\pi_1(M)$.
This group acts on the manifold $M$
with quotient $O''$. Now, $\pi_1(O'')$ contains
${\Bbb Z}/2{\Bbb Z} \times {\Bbb Z}/2{\Bbb Z}$, and hence some singular
point of $O''$ has local group that contains
${\Bbb Z}/2{\Bbb Z} \times {\Bbb Z}/2{\Bbb Z}$. The group of
covering transformations must contain the local group
of this vertex. Hence,
$\pi_1(O'')/\pi_1(M)$ contains ${\Bbb Z}/2{\Bbb Z}
\times {\Bbb Z}/2{\Bbb Z}$. Let $h_1$ and $h_2$ be
the corresponding commuting covering transformations
of $M$, which are involutions.
Let $h_3$ be the composition of $h_1$ and $h_2$,
which also is an involution. 
Let $O_i$ be the quotient $M/h_i$.
Since $h_i$ has non-empty fixed point set,
${\rm sing}(O_i)$ is a non-empty collection of
simple closed curves with order 2.

We claim that, for $i = 1$, 2 or 3, $b_1(O_i) \geq 2$.
Now, each $h_i$ induces an automorphism $h_{i \ast}$
of $H_1(M; {\Bbb R})$. It is clear that
$b_1(O_i)$ is equal to the
dimension of the $+1$ eigenspace of $h_i$.
Suppose that this is at most 1 for
$i = 1$ and $2$. Then the dimension of
the $-1$ eigenspace is at least 3 for
$i = 1$ and $2$. Hence, the intersection of
these eigenspaces has dimension at least 2.
This lies in the $+1$ eigenspace for
$h_{3 \ast}$, and so $b_1(O_3) \geq 2$,
proving the claim.

Hence, by Theorem 3.8 and the Remark preceeding the preceeding this proof,
$\pi_1(O'')$ is large and hence so is $\pi_1(O)$. $\square$

\noindent {\bf Theorem 1.13.} {\sl Let $O$ be a closed orientable hyperbolic
3-orbifold containing
${\Bbb Z}/2{\Bbb Z} \times {\Bbb Z}/2{\Bbb Z}$ in its fundamental group. 
Suppose that $vb_1(|O|) > 0$. Then $\pi_1(O)$ is large.}

\noindent {\sl Proof.} By Proposition 1.4, $O$ is finitely covered
by an orbifold $O'$ such that ${\rm sing}_2^-(O')$ is
non-empty. Since $vb_1(|O|)$ is positive, so is $vb_1(|O'|)$.
Let $|O''|$ be a cover of $|O'|$ such that $b_1(|O''|) >0$,
and let $O''$ be the corresponding cover of $O'$.
Then ${\rm sing}_2^-(O'')$ is non-empty. Hence, by
Theorem 1.9, $\pi_1(O'')$ is large and therefore so
is $\pi_1(O)$. $\square$

The following is a major conjecture in 3-manifold theory.
It has been verified in many cases. For example, Dunfield
and Thurston have checked that it is true for the 10986
hyperbolic 3-manifolds in the `census'.

\noindent {\bf Conjecture 1.14.} {\sl For any closed orientable
3-manifold $M$ with infinite fundamental group, $vb_1(M) > 0$.}

Our main goal in this section is to show that it has
the following slightly surprising consequence.

\noindent {\bf Theorem 1.15.} {\sl If $vb_1(M) > 0$ for all
closed orientable 3-manifolds $M$ with infinite fundamental group, 
then $\pi_1(N)$ is large for any arithmetic 3-manifold $N$,
and hence $vb_1(N)= \infty$.}

This follows from Theorem 1.13 and the following result.

\noindent {\bf Theorem 1.16.} {\sl Let $O$ be a closed orientable 
hyperbolic 3-orbifold
containing ${\Bbb Z}/2{\Bbb Z} \times {\Bbb Z}/2{\Bbb Z}$ in its fundamental group.
Then $O$ is finitely covered by a 3-orbifold $\tilde O$ containing
${\Bbb Z}/2{\Bbb Z} \times {\Bbb Z}/2{\Bbb Z}$ in its fundamental group and for
which $\pi_1(|\tilde O|)$ is infinite.}

The proof will involve an application of the Golod-Shafarevich
theorem, which asserts that groups having a certain type of
presentation are infinite. We therefore need the following lemma.

\noindent {\bf Lemma 1.17.} {\sl Let $O$ be a closed
orientable 3-orbifold with non-empty singular locus. 
Then $\pi_1(O)$ has a finite presentation
$\langle X | R \rangle$ with $|R| - |X| \leq 2
b_1({\rm sing}(O)) - 2$.}

\noindent {\sl Proof.} Let $M$ be the manifold 
obtained from $O$ by removing an open regular neighbourhood
of its singular locus. Then $\pi_1(M)$ has a presentation
$\langle X | R' \rangle$ with $|R'| - |X| {1 \over 2}\chi(\partial M) - 1 = \chi({\rm sing}^-(O)) - 1$. We obtain
$\pi_1(O)$ from $\pi_1(M)$ by adding relations
that are powers of the meridians of the singular locus.
For each circle component of the singular locus,
there is one such meridian. For each graph component $Y$,
the number of meridians is equal to $-3 \chi(Y)$.
Hence, $\pi_1(O)$ has a presentation
$\langle X | R \rangle$ with 
$$|R| - |X| = -2\chi({\rm sing}^-(O)) + |{\rm sing}^0(O)| -1
\leq 2b_1({\rm sing}(O)) - 2.$$
$\square$

\noindent {\bf Proposition 1.18.} {\sl Let $O$ be a closed orientable
hyperbolic 3-orbifold
containing ${\Bbb Z}/2{\Bbb Z} \times {\Bbb Z}/2{\Bbb Z}$ in its fundamental group.
Let $n$ be any integer.
Then $O$ is finitely covered by a 3-orbifold $\tilde O$ for
which ${\rm sing}_2^-(\tilde O)$ is non-empty and
$d_2(\tilde O) \geq n$.}

For this, we will need the following, which is
due to Lubotzky, and appears as Theorem 3.2 in [].

\noindent {\bf Theorem 1.19.} {\sl Let $O$ be an orientable
finite-volume hyperbolic 3-orbifold (with possibly empty
singular locus). Then, for any prime $p$,
$$\sup \{ d_p(K): K \hbox{ is a finite index normal
subgroup of } \pi_1(O)\} = \infty.$$}

\noindent {\sl Proof of Proposition 1.18.} By 
Theorem 1.19, $\pi_1(O)$ has a finite index
normal subgroup $K$ such that $d_2(K) \geq 4n-3$.
We know that $\pi_1(O)$ contains a copy $A$
of ${\Bbb Z}/2{\Bbb Z} \times {\Bbb Z}/2{\Bbb Z}$. Let
$\tilde O$ be the covering space of $O$ corresponding to the
subgroup $KA$. Then, as in the proof of Proposition 1.4,
${\rm sing}_2^-(\tilde O)$ is non-empty.
We claim also that $d_2(KA) \geq n$.
This is because $K$ is a normal subgroup of 
$KA$, and the quotient
is isomorphic to either the trivial group, 
${\Bbb Z}/2{\Bbb Z}$ or ${\Bbb Z}/2{\Bbb Z} \times {\Bbb Z}/2{\Bbb Z}$.
Proposition ? of [] states that
$$d_2(K) - 1 \leq [KA:K] (d_2(KA) - 1),$$
and hence
$$d_2(K A) \geq {d_2(K) - 1 \over 4} + 1\geq n.$$
$\square$

\noindent {\bf Proposition 1.20.} {\sl Let $X$ be a
finite trivalent graph such that $b_1(X) \geq 2$.
Then $X$ contains a subgraph $Y$ such that
$b_1(Y) = 2$ and the number of edges of $Y$ is
at most $6\log_2(b_1(X) - 1) + 8$.}

\noindent {\sl Proof.} We may assume that $X$ is connected. Consider the path
metric on $X$,
where each edge has length 1. For any vertex $v$ and non-negative integer 
$n \leq b_1(X)$, let 
$R_n(v)$ be the minimal radius of a ball centred at $v$
that contains a subgraph $Y$ with $b_1(Y) \geq n$.
Fix a vertex $v$ where $R_2(v)$ has minimal value, $R$ say.
We claim that $$R \leq \log_2(b_1(X) - 1) + 2.$$
There are two cases to consider: when $R_1(v) < R_2(v)$
and when $R_1(v) = R_2(v)$. Let us concentrate on the first case. 

For any non-negative integer $r$, let $B(r)$ be the
ball of radius $r$ around $v$. Now, $B(R_1(v) - 1)$
has trivial $b_1$, and is therefore a tree. But,
by definition of $R_1(v)$, two vertices
in $B(R_1(v) - 1)$ are joined by an edge $e$ in $X - B(R_1(v) - 1)$.
Since $R_1(v) < R_2(v)$, there is a unique such
edge $e$. Let $C$ be the circle in $B(R_1(v))$
consisting of $e$ and the unique embedded path in
$B(R_1(v)-1)$ joining the endpoints of $e$.
Then $C$ consists of at most $2 R_1(v) - 1$ edges.
Now, $B(R - 1)$ is obtained from $C$ by
attaching trees. The distance from any vertex of $B(R - 1)$ to
$C$ is at most $2R - 2$. By definition of $R = R_2(v)$,
there exist two different vertices in $B(R - 1)$
joined by an edge $e'$ in $X - B(R-1)$.
Let $p_1$ be the shortest path in $B(R-1)$ from one of these
vertices to $C$, and let $p_2$ be the shortest
path in $B(R-1)$ from the other vertex to $C \cup p_1$.
Set $Y = e' \cup p_1 \cup p_2 \cup C$.
Then it is clear that $b_1(Y) = 2$.
By construction, $Y$ consists of at most
$6R - 4$ edges.

Now, $C$ runs through at most two of the three
edges adjacent to $v$. Thus, if one were to
remove $v$ and its adjacent edges from $B(R - 1)$, one component would
consist of a based binary tree. The number of vertices in this
tree is equal to $2^{R-1}$. This is a lower bound
for the number of vertices in $X$, which we denote by $|V(X)|$.
But, since $X$ is trivalent, $b_1(X) = {1 \over 2}|V(X)| + 1$.
So, $$2^{R - 1} \leq |V(X)| = 2 b_1(X) - 2$$
and therefore
$$R \leq \log_2(b_1(X) - 1) + 2.$$
This proves the claim when $R_1(v) < R_2(v)$.
When these quantities are equal, the proof is similar
but simpler, and so is omitted.

We have already seen that $Y$ contains at most
$6R - 4$ edges, and so the proposition is proved.
$\square$

The final step in the proof of Theorem 1.16 is the
following proposition.

\noindent {\bf Proposition 1.21.} {\sl Let $O$ be a 
closed orientable 3-orbifold with fundamental group $G$. 
Let $Y$ be a subgraph of ${\rm sing}_2(O)$ that has
$b_1(Y) = 2$ and at most $6\log_2(b_1({\rm sing}_2(O)) - 1) + 8$
edges. Let $K$ be the subgroup of $G$ normally generated
by the meridians that encircle $Y$. Let $\tilde O$ be the
covering space of $O$ corresponding to the subgroup $K G^2$. 
Suppose that $d_2(O) \geq N$. [{\bf $N$ is an integer yet to
be calculated.}]
Then $\pi_1(|\tilde O|)$ is infinite.}

\noindent {\sl Proof.} Note first that it
suffices to consider orbifolds where every arc
and circle of the singular locus has singularity
order 2. For, we may let
$O_1$ be the orbifold obtained from $O$ by removing
each component of the singular locus with odd
order and replacing it with manifold points,
and relabelling each even singularity order by 2.
Let $G_1 = \pi_1(O_1)$.
There is an obvious map $q \colon O \rightarrow O_1$
that induces a surjection $q_\ast \colon \pi_1(O) \rightarrow \pi_1(O_1)$
and an isomorphism $H_1(O; {\Bbb Z}/2{\Bbb Z}) \rightarrow 
H_1(O_1; {\Bbb Z}/2{\Bbb Z})$. Let $Y_1$ be the image of $Y$
under $q$, and let $K_1$ be the subgroup of 
$G_1$ normally generated by the meridians of $Y_1$.
Then $K_1 G_1^2$ is the image of $K G^2$ under $q_\ast$. 
Let $\tilde O_1$ be the covering space of $O_1$
corresponding to $K_1 G_1^2$. Then, $\tilde O_1$
and $\tilde O$ have the same underlying manifold.
Therefore, if we can show that $\pi_1(|\tilde O_1|)$ is
infinite, then so is $\pi_1(|\tilde O|)$.

Let $M_1$ be the meridians of $Y$. 
Let $M_2$ be those meridians
of ${\rm sing}(O)$ that lie in $G^2$
but that are not meridians of $Y$. 
Then $|M_1| + |M_2|$ is at most the total
number of meridians of $Y$, which is at most
${3 \over 2}b_1({\rm sing}(O))$.
Consider the group $\Gamma\pi_1(O)/\langle \! \langle M_1, M_2 \rangle \! \rangle$.
Lemma 1.17 states that $\pi_1(O)$ has a finite presentation
where the number of relations minus the number of generators is
at most $2b_1({\rm sing}(O)) -2$.
Hence, $\Gamma$ has a presentation $\langle X | R \rangle$
with $|R| - |X| \leq 2 b_1({\rm sing}(O)) - 2 + |M_1| + |M_2|
\leq {7 \over 2} b_1({\rm sing}(O)) - 2$. By Proposition 1.1,
this is at most ${7 \over 2} d_2(O) - 2$. Now, adding the
relations in $M_2$ does not affect $d_2$, because
they lie in $G^2$. Hence, 
$$d_2(O) \geq d_2(\Gamma) \geq d_2(O) - |M_1| \geq d_2(O) - 6\log_2(d_2(O) - 1)
- 8.$$
Therefore,
$$\eqalign{
&d_2(\Gamma)^2/4 - |R| + |X| - d_2(\Gamma)\cr
&\qquad \geq d_2(\Gamma)^2/4 - (7/2) d_2(O) + 2 - d_2(\Gamma) \cr
&\qquad \geq (d_2(O) - 6\log_2(d_2(O) - 1) - 8)^2/4 - (9/2) d_2(O) + 2\cr
&\qquad > 0,}$$
by our hypothesis that $d_2(O) > N$.
Hence, by the Golod-Shafarevich theorem, $\Gamma$ is
infinite, and hence so is $\Gamma^2$, because it has
finite index in $\Gamma$.

Let $\tilde O$ be the covering space of $O$ corresponding
to $K G^2$. We claim that $\pi_1(|\tilde O|) \cong \Gamma^2$.
This will prove that $\pi_1(|\tilde O|)$ is infinite.
Now, $\pi_1(|\tilde O|)$ is obtained from $\pi_1(\tilde O)
= K G^2$ by quotienting each meridian of the singular
locus. But the meridians in $O$ that lift to
meridians of singular components of $\tilde O$ are
precisely $M_1 \cup M_2$. Hence, $\pi_1(|\tilde O|)
= K G^2 / \langle \! \langle M_1, M_2 \rangle \! \rangle
= \Gamma^2$. $\square$

\noindent {\sl Proof of Theorem 1.16.} 
By Proposition 1.18, $O$ is finitely covered by
an orbifold $O'$ for which ${\rm sing}_2^-(O')$ is
non-empty and for which $d_2(O') \geq N$.
By Proposition 1.20, ${\rm sing}_2^-(O')$ contains
a subgraph $Y$ with $b_1(Y) = 2$ and with at
most $6\log_2(b_1({\rm sing}_2^-(O')) - 1) + 8$ edges,
which is at most $6\log_2(2 b_1({\rm sing}_2(O')) - 1) + 8$.
Let $O''$ be the covering space of $O'$
corresponding to $K G^2$, where $G = \pi_1(O')$
and $K$ is the subgroup normally generated
by the meridians of $Y$. Since it contains $K$,
the inverse image of $Y$ in $O''$ is a disjoint
union of copies of $Y$. In particular,
${\rm sing}_2^-(O'')$ is non-empty, and so $\pi_1(O'')$
contains ${\Bbb Z}/2{\Bbb Z} \times {\Bbb Z}/2{\Bbb Z}$. By
Proposition 1.21, $\pi_1(|O''|)$ is infinite.
Thus, by our hypothesis, $vb_1(|O''|) > 0$.
So, applying Theorem 1.13, $\pi_1(O'')$ is large
and hence so is $\pi_1(O)$. $\square$

\noindent {\bf Conjecture.} (Lubotzky-Sarnak) {\sl
The fundamental group of any closed hyperbolic 3-manifold
does not have property $(\tau)$.}

I can in fact prove:

\noindent {\bf Theorem.} {\sl The Geometrisation Conjecture
and Lubotzky-Sarnak Conjecture imply that the fundamental
group of any arithmetic 3-manifold is large.}

The proof requires a new version of Theorem 1.5 specifically
tailored to 3-orbifolds. It would take about a couple of pages to write
down, and so might not be worth it?

that $\Gamma$ contains $A_4$. 
If $\nu$ in $\Ram_f(B)$, then $k_\nu$ contains no quadratic extension of ${\Bbb Q}_2$.}

Note that $A_4$ is a subgroup of $S_4$ and $A_5$. The
following is contained in [MR] (see Theorem 5.4.1).

\noindent {\bf Theorem 2.4}~{\it Suppose that $\Gamma$ is a Kleinian group of finite co-volume 
with invariant trace-field and quaternion algebra denoted by $k$ and $B$ respectively. Assume
that $\Gamma$ contains $A_4$. Then,
$$B \cong \biggl ({{-1,-1}\over k}\biggr).$$
Furthermore, if $\Gamma$ contains $A_5$, then ${\Bbb Q}(\sqrt{5}) \subset k$.}

Note that the only finite places that can ramify $B$ (as in Theorem 2.4) are dyadic.  

In the arithmetic setting, one can say more, and which will be of use below.

3-orbifolds.

\noindent {\bf Theorem 6.6.} {\sl Let $\Gamma$ be an arithmetic
Kleinian group, commensurable with a group derived from a quaternion algebra
containing a dihedral subgroup. Then $\Gamma$ is large.}

\noindent {\sl Proof.}

Only elementary finite group theory is required here, together with
the fact that 
$$\{ A \in {\rm GL}(2, \hat {\Bbb Z}_p): A \equiv I \hbox{ mod } p^i \}$$
is a uniform pro-$p$ group for each $i \geq 2$, and $i = 1$ when $p$ is odd
(where $\hat {\Bbb Z}_p$ is the $p$-adic integers).
This is a relatively easy result in the theory of pro-$p$ groups
(Theorem 7.2 in []).

At one point in the argument, use is made of the
Prime Number Theorem. But, that aside, the deduction
of Theorem 9.3 is a short application of fairly elementary
group theory. In particular, the Golod-Shafarevich inequality
and the full theory of $p$-adic analytic groups and  are not required.

There is a slight variation on the above proof of Theorem 9.1,
which we merely sketch here. Instead of counting all finite index
subgroups of $\Gamma$, one can focus on those subgroups
of $\Gamma$ which are subnormal. Therefore, let $s_n^{\triangleleft \triangleleft}(\Gamma)$
and $c_n^{\triangleleft \triangleleft}(\Gamma)$ be the number of
subnormal subgroups (respectively, subnormal congruence subgroups)
of $\Gamma$ with index at most $n$. By slightly extending
the methods of Section 2, one can show that $s_n^{\triangleleft \triangleleft}(\Gamma)$
grows exponentially as a function of $n$. (See Theorem 1.4 in []).
Thus, we need only show that $c_n^{\triangleleft \triangleleft}(\Gamma)$
grows more slowly than this. Because subnormal subgroups of
${\cal O}^1/{\cal O}^1(m{\cal O})$ are easier to count than
arbitrary subgroups, some of the technicalities in the proof of Theorem 9.3
(including the use of the Prime Number Theorem) can thereby be avoided.